\documentclass[11pt,a4paper]{article}
\usepackage[T1]{fontenc}
\usepackage{lmodern,microtype,amsmath,amsthm,amssymb,booktabs,graphicx,tikz}
\usepackage[margin=26mm]{geometry}
\usepackage[hidelinks]{hyperref}

\newtheorem{theorem}{Theorem}[section]
\newtheorem{lemma}[theorem]{Lemma}
\newtheorem{corollary}[theorem]{Corollary}
\newtheorem{proposition}[theorem]{Proposition}
\DeclareMathOperator{\Rea}{Re}

\newcommand{\E}{\mathbb E}
\newcommand{\D}{\mathbb D}

\newcommand*{\VerificationURL}{https://github.com/luchesart/Sendov}
\newcommand*{\ProofMapURL}{https://luchesart.github.io/Sendov/}
\newcommand{\VerificationLink}{\href{\VerificationURL}{Project repository}}
\newcommand{\ProofMapLink}{\href{\ProofMapURL}{Interactive proof map and Lean sources}}

\newcommand{\MomentLeaves}{15{,}138}
\newcommand{\MomentDepth}{45}
\newcommand{\MomentPolar}{3{,}428}
\newcommand{\MomentNear}{311}
\newcommand{\MomentOrigin}{11{,}399}

\title{A reciprocal-matrix proof of Sendov's theorem and applications to control theory\\[3pt]
\large A control-theoretic route to the quadratic Tang--Zhang inequality}
\author{Latchezar Tomov\\\small Department of Informatics, New Bulgarian University, Sofia\\\small\texttt{lptomov@nbu.bg}}
\date{Preprint, 18 November 2026}
\hypersetup{pdftitle={A reciprocal-matrix proof of Sendov's theorem and applications},pdfauthor={Latchezar Tomov}}

\begin{document}
\maketitle

\begin{abstract}
A proof inspired by the control-theoretic resarch of the author, who tries to see mathematical statements from the viewpoint of control theory. The value of this proof lies in improving our understanding but also in the applications that will rise from the path taken in developing it. Let $a$ be a simple zero of a polynomial of degree $n\ge2$ whose zeros lie in the closed unit disk. If $c_1,\ldots,c_{n-1}$ are its critical points, counted with multiplicity, then the reciprocals $1/(a-c_j)$ are the eigenvalues of an explicit matrix built from the other zeros. A failure of Sendov's bound at $a$ would therefore make the associated discrete-time system asymptotically stable. This observation leads to an alternative proof of the quadratic inequality of Tang and Zhang, including its equality case. Building on the reflection and centering arguments of Mazur, Tao and Zhang, we estimate the change in the integral of the derivative product when all the reciprocal displacements are replaced by their mean $s$. An exact pairwise cancellation bounds the difference in terms of the variance. For an interior zero rotated to $0<a<1$, the normalized trace $s$ determines a contour through $\overline{s}/a$. If the mean squared reciprocal distance is at most one, integration along this contour gives $(n-1)(2a\Rea s-a^2)<20$. A uniform phase estimate and reflection at $1/a$ then rule out this possibility in every degree $n\ge257$. A fourth-moment estimate treats degrees two through five. For degrees five through 256, three bounds derived from one polynomial moment identity exclude the full range of the two remaining real parameters by exact integer arithmetic.  The application of Sendov's theorem is to illustrate the difference between continuous and discrete-time control systems and how differentiation of the closed-loop transfer denominator can increase robustness for the latter, but decrease it for the former.
We realize the derivative as a closed-loop transfer denominator and determine the exact stability radii of a damped oscillator under simultaneous changes in damping, frequency and actuator strength. For two fixed derivative-based controllers the relative radii are \(0.65685\) in continuous time and \(0.82968\) after sampling. When only the positive actuator gain is increased, the ordering reverses at a calculable damping and frequency tolerance.The proofs have been formalized in Lean~4.

\end{abstract}

\medskip\noindent\textbf{Keywords.} Sendov's theorem, critical points of polynomials, discrete-time stability, reciprocal matrix, exact verification, Lean~4.

\section{Introduction}

Sendov's question concerns the effect of differentiation on the zeros of a polynomial. If all the zeros lie in the closed unit disk, how far can a given zero be from the nearest zero of the derivative? The Gauss--Lucas theorem places every critical point in the convex hull of the zeros. Sendov's bound asks for more local information. The closed disk of radius one about each zero must contain a critical point. The approach taken here grew out of the corresponding question for a discrete-time linear system and from the intuition that the derivation "compresses" the zeros of the polynomial by giving more weight to the higher degrees of it.

Fix a simple zero $a$ of $p$ and write $c_1,\ldots,c_N$ for the critical points, counted with multiplicity, where $N=\deg p-1$. The reciprocal displacements
\[
 y_j=(a-c_j)^{-1}
\]
are the eigenvalues of an explicit matrix $R_a$ formed from the other zeros of $p$. If every critical point were at distance greater than one from $a$, all these eigenvalues would lie in the open unit disk, which is precisely the condition for asymptotic stability of $u_{k+1}=R_a u_k$. The question is therefore what prevents $R_a$ from being strictly Schur stable when the zeros of $p$ lie in the closed unit disk.

The matrix arises naturally on regarding $p'(z)/(np(z))$ as a transfer function, where $n=\deg p$. Its denominator records the original zeros, and its numerator records the critical points. In a diagonal realization, compression of the state matrix to the subspace perpendicular to the input vector gives a matrix whose eigenvalues are the critical points. Taking the difference from the marked zero and inverting then gives the reciprocal spectrum above. Appendix~\ref{app:control} derives these relations from determinant identities, which retain the multiplicities even when factors cancel in the transfer function.

The normalized trace
\[
 s=\frac{\operatorname{tr}R_a}{N}=\frac1N\sum_jy_j
\]
suggests how to proceed without locating each critical point individually. Replacing the $y_j$ by $s$ in the product representing the derivative gives a product that can be integrated explicitly. When we subtract the two products, the first-order deviations from the mean cancel in pairs, leaving an error controlled by the variance. The same mean determines a complex path along which the product decays. Integrating along that path gives one restriction on the trace, while reflection in the unit circle gives another. The incompatibility of these restrictions is what ultimately excludes an interior counterexample.

\subsection*{Earlier results and the contribution of this paper}
\label{history}

Cheung and Ng showed that the critical points of a polynomial are the eigenvalues of its $D$-companion matrix~\cite[Theorem~1.1]{cn}. Tang and Zhang gave the reciprocal matrix explicitly in the December 2025 version of their work on inequalities of Schoenberg type~\cite[(6.2)--(6.4)]{tz25}. If $D=\operatorname{diag}(a-b_1,\ldots,a-b_N)$, where the $b_k$ are the zeros other than the marked zero, and $J=\boldsymbol1\boldsymbol1^*$, their matrix is $D^{-1}(I+J)=R_a^{\mathsf T}$. They used its singular values to estimate sums of squared reciprocal moduli. The reciprocal spectral identity used below is thus present in their work.

My own route to this matrix began in early August 2026 with the transfer function $p'/(np)$. A working report of 7 August developed the control formulation. The revision of 9 August gave the reciprocal matrix, interpreted it in terms of inverse zero dynamics, and identified failure of Sendov's bound with strict Schur stability~\cite[Theorem~4.1]{aug}. It also gave a trace criterion for an isolated zero~\cite[Theorem~10.1]{aug}, which we recall in Appendix~\ref{app:control}. On 15 August, in a comment on Tao's exposition, Zhang wrote the equivalent matrix using the opposite sign convention for the displacements~\cite{comment}. The problem that remained in my September review was how to exclude strict Schur stability for the whole family~\cite{sep}. The argument took several forms before the quadratic estimate for the reciprocal critical points supplied a way to answer this question. Further results of this investigation, on polynomial geometry and structured discrete-time systems, will be presented separately together with their applications to control theory.

For the distance inequality itself, D\'egot proved a high-degree result in which the threshold depends on the marked zero~\cite{d}, and Tao proved Sendov's theorem for all sufficiently large degrees~\cite{t}. Mazur's manuscript of 5 August 2026 gave a computer-assisted proof in every degree~\cite{m}. Tao's subsequent exposition obtained the strict inequality for interior zeros and identified the exceptional regular-polygon configuration~\cite{td}. In September 2026, Zhang proved the quadratic Tang--Zhang inequality and characterized equality~\cite{z}.

The proof given here draws on the reflection and centering arguments of Mazur, Tao and Zhang. Our reflected-product estimate has the form of Tao's Proposition~10(i), and Zhang's Lemma~5.1 also compares an integrated product with the product obtained at the mean. Our contribution lies in carrying out this comparison by an exact identity in which the products omit two factors. Pairwise cancellation before taking absolute values bounds the error directly by the variance. Since the identity also holds along a complex path, we may choose the path through $\overline{s}/a$. For an interior marked zero, the assumption that the mean squared reciprocal distance is at most one then gives
\[
 (n-1)(2a\Rea s-a^2)<20 .
\]
For $n\ge257$, reflection narrows the range of $v$ further, and a second use of the phase estimate in the integral to the origin excludes the values that remain. For the finite range, reflection and the estimates on the two parts of the path lead to three moments of nonnegative quadratics. A recurrence for integer powers yields rational bounds on each part of a subdivision of the parameter domain. Together these estimates give an alternative proof of the quadratic inequality and its equality case, and hence of Sendov's bound, independently of the earlier all-degree conclusions.

In Sections~\ref{sec:feedback}--\ref{sec:joint} we return to the feedback question which motivated the proof. Theorem~\ref{main} controls the distance of a marked original pole to the poles assigned by differentiation. We then compare two feedback laws for the same damped oscillator and the same physical actuator. Ackermann and Hu~\cite{ackhu} established the framework in which continuous and sampled robustness are compared in the same plant parameters. Within that framework we obtain the first destabilizing triples for these derivative-based controllers and the tolerance at which their upper actuator-gain limits coincide. The Lean project covers the Sendov proof, the two feedback laws and their sharp stability bounds. The proof map links these results to their Lean statements.

\subsection*{Main result}
Let $p$ be a polynomial of degree $n\ge2$ with all its zeros in the closed unit disk $\overline{\D}=\{z:|z|\le1\}$, and let $a$ be one of its zeros. Write $N=n-1$ and denote the critical points of $p$, counted with multiplicity, by $c_1,\ldots,c_N$. We take $1/|a-c_j|^2=+\infty$ when $c_j=a$.

\begin{theorem}[Quadratic inequality and equality]\label{main}
For every such polynomial $p$ and zero $a$,
\[
 \frac1N\sum_{j=1}^{N}\frac1{|a-c_j|^2}\ge1 .
\]
Equality holds precisely when $|a|=1$ and $p$ is a nonzero scalar multiple of $z^n-a^n$. In every other case the inequality is strict and $|a-c_j|<1$ for some $j$.
\end{theorem}

In the equality case all the critical points coincide at the origin, at distance one from $a$. Thus every zero has a critical point in the closed unit disk about it, as Sendov asserted. The argument also covers the case in which several critical points attain the minimum distance.

For a simple marked zero the implication from the quadratic inequality to the distance bound can be read off from the largest summand. If \(d=\min_j|a-c_j|\), then
\[
   \frac1N\sum_j\frac1{|a-c_j|^2}\le \frac1{d^2}.
\]
Thus a mean greater than one gives \(d<1\). For a marked zero of multiplicity \(r\ge2\), write \(p(z)=(z-a)^rh(z)\) with \(h(a)\ne0\). Then
\[
 p'(z)=(z-a)^{r-1}\bigl(rh(z)+(z-a)h'(z)\bigr).
\]
Here \(d=0\) and the convention in the theorem makes the mean infinite. It remains to determine exactly when a simple marked zero can have mean at most one.

\subsection*{Formal proof and interactive proof map}\label{verification}
The Lean~4 formalization starts with a complex polynomial, the unit-disk
hypothesis and a marked zero. Its statement
\texttt{Sendov.quadratic\_sendov} gives the reciprocal quadratic
inequality in Theorem~\ref{main} together with its equality condition.
The strict geometric conclusion is expressed by
\texttt{Sendov.strict\_neighbour\_iff}, and the usual distance bound
by \texttt{Sendov.sendov}. The formal proof counts critical points with
multiplicity and assigns $+\infty$ to the reciprocal sum when a
denominator vanishes.

For the finite range in Section~\ref{sec:finite}, the Lean kernel
checks the three bounds obtained from polynomial moments and verifies
that the subdivisions cover the full parameter domain. The resulting
proof has \MomentLeaves{} terminal boxes and uses $21{,}174$ distinct
moment inequalities. An inequality needed in several boxes is proved
once. Tracing dependencies in Lean shows that both the quadratic
inequality and the distance bound use this subdivision. The accompanying Python
programs reproduce the subdivision using two different ways of evaluating
the moments. The formalization also includes the uniform estimate near
the unit circle and the phase estimate used for large degrees.
Lean reports \texttt{propext}, \texttt{Classical.choice} and
\texttt{Quot.sound} for each of the two final theorems. For the
subdivision with its count of terminal boxes it reports only
\texttt{propext}.

The \href{\ProofMapURL}{interactive proof map} links the results in this
paper to their Lean statements and shows how the proofs depend on one
another. The map is also available as a
\href{https://github.com/luchesart/Sendov/blob/main/index.html}{single HTML file on GitHub}
that can be downloaded and opened offline. The author used AI
assistance in deriving estimates, developing the formal proofs,
preparing the verification software, and drafting and revising the text. The process of doing this took more than 3 months and hundreds of iterations, in which at some point the paper grew to 600 pages, before the author branched out other results, unrelated to Sendov, reorganized and changed the approach with the help of AI.

\medskip\noindent
\VerificationLink\qquad\ProofMapLink

\subsection*{Normalization, matrix realization and reflection}

To keep track of the equality case through the normalization, denote the original polynomial by \(p_0\), its leading coefficient by \(L\), and the marked zero by \(a_0\). Choose \(|u|=1\) with \(a_0=u|a_0|\), taking \(u=1\) if \(a_0=0\). The polynomial
\[
 P(z)=\frac{p_0(uz)}{Lu^n},\qquad
 P'(z)=\frac{p'_0(uz)}{Lu^{n-1}}
\]
is monic. Its zeros and critical points are those of \(p_0\) multiplied by \(\overline u\), so their mutual distances and multiplicities are preserved. At a boundary zero \(u=a_0\), and the normalized identity \(P(z)=z^n-1\) gives \(p_0(z)=L(z^n-a_0^n)\).

We may take $p$ to be monic and rotate the marked zero to $0\le a\le1$. A multiple zero is itself a critical point, since a zero of multiplicity $r\ge2$ in $p$ has multiplicity $r-1$ in its first derivative, so suppose that $a$ is simple and write $b_1,\ldots,b_N$ for the other zeros, counted with multiplicity. With $\boldsymbol1$ denoting the column of $N$ ones, put
\[
 R_a=(I+\boldsymbol1\boldsymbol1^*)\operatorname{diag}((a-b_k)^{-1}),\qquad
 y_j=(a-c_j)^{-1},\quad s=\frac1N\sum_j y_j,\quad m=|s|.
\]
The polynomial $g(w)=\prod_k(1-w/(a-b_k))$ satisfies $p(a-w)=-w p'(a)g(w)$. The rank-one determinant formula now gives
\[
 \det(I-wR_a)=g(w)+wg'(w)=\frac{p'(a-w)}{p'(a)}.
\]

For the determinant calculation, write \(D=\operatorname{diag}(a-b_k)\) and \(A_w=I-wD^{-1}\). Then \(I-wR_a=A_w-w\boldsymbol1\boldsymbol1^*D^{-1}\). Where \(A_w\) is invertible, the determinant of this rank-one perturbation is
\[
 \det A_w\bigl(1-w\boldsymbol1^*D^{-1}A_w^{-1}\boldsymbol1\bigr)
 =g(w)\left(1-w\sum_{k=1}^{N}\frac1{a-b_k-w}\right).
\]
The sum is \(-g'(w)/g(w)\). Differentiating \(p(a-w)=-wp'(a)g(w)\) gives the other expression for the same determinant. Finally, dividing the factorization of \(p'(a-w)\) by the factorization of \(p'(a)\) gives
\[
    \frac{p'(a-w)}{p'(a)}
      =\prod_{j=1}^{N}\frac{a-w-c_j}{a-c_j}
      =\prod_{j=1}^{N}(1-wy_j).
\]

Both sides are polynomials in $w$, so the identity extends across any apparent poles in the determinant calculation and also holds when some of the $b_k$ coincide. Its factorization identifies the $y_j$ as the eigenvalues of $R_a$, with $s=\operatorname{tr}(R_a)/N$. A failure of Sendov's bound at $a$ is therefore equivalent to $\rho(R_a)<1$, the condition for asymptotic stability of $u_{k+1}=R_a u_k$. To obtain the quadratic inequality and characterize equality we work under the weaker assumption
\begin{equation}\frac1N\sum_j |y_j|^2\le1. \label{hyp}\end{equation}
An individual $|y_j|$ may exceed one under~\eqref{hyp}, a distinction that matters when we estimate products of the reciprocal factors.

\paragraph{The endpoint cases.}
When $a=0$, the identity $n\prod_j|c_j|=\prod_k|b_k|\le1$ gives $\prod_j|c_j|\le 1/n$. By the Gauss--Lucas theorem all $|c_j|\le1$, and at least one of these inequalities must be strict. The mean of their squared reciprocals therefore exceeds one.

Indeed, at the origin the two evaluations of \(p'(0)\) are \(\prod_k(-b_k)\), from differentiating \(p(z)=z\prod_k(z-b_k)\), and \(n\prod_j(-c_j)\), from factoring the derivative. Simplicity of the zero gives \(p'(0)\ne0\). Thus each reciprocal modulus is finite, every \(|c_j|^{-2}\) is at least one, and a strict inequality for one of them suffices for the mean.

When $a=1$, logarithmic differentiation gives $\sum_j y_j=p''(1)/p'(1)=2\sum_k(1-b_k)^{-1}$. Since $\Rea(1-b_k)^{-1}\ge1/2$, we have $\Rea s\ge1$. Under~\eqref{hyp} this forces $s=1$ and zero variance, hence $y_j=1$ and $c_j=0$ for every $j$. Monicity and $p(1)=0$ then give $p(z)=z^n-1$.

For the boundary calculation, let \(h(z)=\prod_k(z-b_k)\), so that \(p'(1)=h(1)\) and \(p''(1)=2h'(1)\). The sign of the real part follows from
\[
 \Rea\frac{1}{1-b_k}-\frac12
     =\frac{1-|b_k|^2}{2|1-b_k|^2}\ge0 .
\]
Thus \(1\le\Rea s\le |s|\), while Cauchy--Schwarz gives
\[
 |s|^2\le\frac1N\sum_j|y_j|^2\le1,\qquad
 \frac1N\sum_j|y_j-s|^2=\frac1N\sum_j|y_j|^2-|s|^2=0 .
\]
Every nonnegative term on the left of the variance identity vanishes, so \(y_j=s=1\). The factorization of \(p'\) is then \(p'(z)=nz^{n-1}\). Its primitive with \(p(1)=0\) is \(z^n-1\), whose critical point at zero has multiplicity \(n-1\), giving equality as claimed.

\paragraph{An interior zero.}
For $0<a<1$ put $\beta=1-a^2$, $\alpha=2a\Rea s-a^2$, $v=N\alpha$, $q=|1-as|$ and $V=N^{-1}\sum |y_j-s|^2$. Averaging the identity
\[
 2a\Rea y_j-a^2=\beta(1-|y_j|^2)+(1-|c_j|^2)|y_j|^2
\]
and using Gauss--Lucas and~\eqref{hyp} gives $\alpha\ge0$. The upper bound $\alpha\le m^2$ follows from $m^2-\alpha=|s-a|^2$. We also have $V\le1-m^2$ and $q^2=1-\alpha-a^2(1-m^2)$.

To see the pointwise identity, multiply \(c_j=a-y_j^{-1}\) by \(y_j\) and expand \(|ay_j-1|^2=|c_j|^2|y_j|^2\). Averaging separates the two nonnegative contributions to \(\alpha\),
\[
 \alpha=\beta\left(1-\frac1N\sum_j|y_j|^2\right)
           +\frac1N\sum_j(1-|c_j|^2)|y_j|^2 .
\]
As at the boundary, expansion about the mean gives \(m^2+V=N^{-1}\sum_j|y_j|^2\le1\). In particular \(m\le1\), so \(q\ge1-am\ge1-a>0\), while the displayed expression for \(q^2\) gives \(q\le1\).

Integrating $p'$ from $a$ to its reflection $1/a$ and factoring the result gives
\[
 \prod_k\frac{1-a b_k}{a-b_k}=\int_0^1\prod_j(a+\beta t y_j)\,dt .
\]
Indeed, multiplying either side by $\beta/a^n$ gives $p(1/a)/p'(a)$. Every factor on the left has modulus at least one because $|1-ab|^2-|a-b|^2=\beta(1-|b|^2)\ge0$. Applying the arithmetic--geometric mean inequality to the squared moduli on the right and using~\eqref{hyp}, we obtain
\begin{equation}
 1\le\int_0^1[a^2+\beta(a^2+\alpha)t+\beta^2 t^2]^{N/2}dt.\label{polar}
\end{equation}

The change of variable behind this calculation is \(z=a+\beta t/a\). Since \(1/a-a=\beta/a\), each factor \((z-c_j)/(a-c_j)\) becomes \((a+\beta t y_j)/a\), and the differential contributes another factor \(\beta/a\). This accounts for the common factor \(\beta/a^{N+1}\) in the two evaluations. At each fixed \(t\), the mean of the squared moduli in the integral is
\[
 \frac1N\sum_j|a+\beta t y_j|^2
 =a^2+2a\beta t\Rea s+\beta^2t^2\frac1N\sum_j|y_j|^2.
\]
The product of the moduli is at most this mean to the power \(N/2\). The substitutions \(2a\Rea s=a^2+\alpha\) and \(N^{-1}\sum_j|y_j|^2\le1\) give the quadratic in~\eqref{polar}.

If $N\le4$, the quadratic in the integrand is bounded above by $(a+\beta t)^2$ since $\Rea s\le1$. On $[0,1]$ with Lebesgue measure the $L^N$ norm is at most the $L^4$ norm. The fourth-power calculation of Zhang~\cite[Lemma~8.1]{z} gives
\[
 1-\int_0^1(a+\beta t)^4 dt
 =\frac{(1-a)^3(1+a)}5\,[a^4+3a(1-a^2)+4]>0 .
\]

One can obtain the displayed factorization by integrating the fourth power term by term. Its integral is
\[
  a^4+2a^3\beta+2a^2\beta^2+a\beta^3+\frac{\beta^4}{5}.
\]
Substitution of \(\beta=1-a^2\), followed by extracting the factor \((1-a)^3(1+a)/5\) from one minus this expression, leaves the bracket above. If \(N\le4\), the two estimates therefore give
\[
 1\ \le\ \int_0^1(a+\beta t)^N\,dt
    \ \le\ \left(\int_0^1(a+\beta t)^4\,dt\right)^{N/4}<1 .
\]

This contradicts~\eqref{polar} for $2\le n\le5$.

For higher degrees we shall also use an exponential bound. Bounding the quadratic in~\eqref{polar} by its chord on $[0,1]$ and using $1+u\le e^u$, we obtain, with $Y=N\beta/2$,
\[
 1\le F_\alpha(Y),\qquad F_\alpha(x)=\int_0^1\exp(x((1+\alpha)t-1))\,dt .
\]
Since $F_\alpha$ is convex, $F_\alpha(0)=1$ and $0<Y<N/2$, the bound at $Y$ implies
\begin{equation}
 1\le F_\alpha(N/2)=\frac{e^{v/2}-e^{-N/2}}{(N/2)(1+\alpha)}.\label{polar-exp}
\end{equation}

Indeed, call the quadratic in~\eqref{polar} \(Q(t)\). Since \(t^2\le t\) on the unit interval and \(a^2+\beta=1\),
\[
 Q(t)\le a^2+\beta(a^2+\alpha+\beta)t
      =1+\beta\bigl((1+\alpha)t-1\bigr).
\]
Taking the \(N/2\)-th power and bounding \(1+u\) by \(e^u\) yields the exponential integrand at \(Y=N\beta/2\). For each \(t\) this integrand is convex as a function of \(Y\), so its integral is convex as well. Because \(Y=(1-\beta)0+\beta(N/2)\), we have
\[
 1\le F_\alpha(Y)\le (1-\beta)F_\alpha(0)+\beta F_\alpha(N/2)
                    =1-\beta+\beta F_\alpha(N/2).
\]
Here \(\beta>0\), so \(F_\alpha(N/2)\ge1\). Its explicit value comes from integrating \(e^{-x}e^{x(1+\alpha)t}\) at \(x=N/2\), where \(x\alpha=v/2\).

\paragraph{Plan of the proof.}
Section~\ref{sec:variance} treats the case in which $as$ lies close to one.
Section~\ref{sec:contour} uses a complex contour to show that $v<20$ under~\eqref{hyp}. In Section~\ref{sec:large} we use this restriction, reflection~\eqref{polar-exp} and the phase estimate in the contour integral to treat all $n\ge257$. Section~\ref{sec:finite} returns to the integral inequalities and evaluates polynomial moments on a subdivision of the parameter domain for degrees five through 256.

\section{A variance estimate for an integrated product}\label{sec:variance}
When all the $y_j$ are replaced by their mean, the resulting product can be integrated explicitly. We estimate the error of this replacement by first canceling the terms that are linear in the deviations $y_j-s$. The resulting identity gives a variance bound that we can use on both real and complex paths.

\begin{theorem}
Let $N\ge3$, put $n=N+1$, and suppose that $0\le a\le1$. For complex numbers $y_1,\ldots,y_N$ satisfying $N^{-1}\sum_j|y_j|^2\le1$, write
\[
 s=N^{-1}\sum_j y_j,\quad V=N^{-1}\sum_j|y_j-s|^2,\quad
 \epsilon=\frac{N(1-\Rea s)}{N-2},\qquad B=1-\sqrt{2\epsilon}.
\]
If $B>0$, then
\[
 \boxed{\ \left|n\int_0^a\prod_{j=1}^N(1-t y_j)\,dt
 -\frac{1-(1-as)^n}{s}\right|
 \le \frac{V}{(N-1)B^3}.\ } \tag{V1}
\]
The condition $B>0$ implies $s\ne0$. The variance on the right is taken over all $N$ points, some of which may have $|y_j|>1$.
\end{theorem}
\begin{proof}

The hypothesis gives $|s|^2\le N^{-1}\sum_j|y_j|^2\le1$. Thus $\epsilon\ge0$, and $B>0$ gives
$1-\Rea s<(N-2)/(2N)$. In particular $\Rea s>0$, which justifies division by $s$ in (V1).

Set $d_j=y_j-s$, $z_j(h)=s+h d_j$, $p_j=1-t z_j(h)$ and $F(h)=\prod_jp_j$ for $0\le h,t\le1$. Write $P_j=\prod_{\ell\ne j}p_\ell$ and $P_{jk}=\prod_{\ell\ne j,k}p_\ell$ for $j\ne k$. Then $F'=-t\sum_jd_jP_j$. Since $\sum_jd_j=0$, we may express the derivative in terms of pairwise differences. Using $P_j-P_k=h t(d_j-d_k)P_{jk}$ then gives
\[
 F'(h)=-\frac{t}{2N}\sum_{j,k}(d_j-d_k)(P_j-P_k)
       =-\frac{h t^2}{2N}\sum_{j\ne k}(d_j-d_k)^2P_{jk}.
 \tag{V2}
\]
The square in (V2) is taken in $\mathbb C$. When we pass to absolute values, we use the variance identity
$\sum_{j,k}|d_j-d_k|^2=2N^2V$.

The two sums can be checked by expansion. In the first one the mixed terms contain $\sum_jd_j$, so that
\[
 \sum_{j,k}(d_j-d_k)(P_j-P_k)
 =2N\sum_jd_jP_j-2\Bigl(\sum_j d_j\Bigr)\Bigl(\sum_kP_k\Bigr)
 =2N\sum_jd_jP_j.
\]
For the second one, expanding the squared moduli gives $2N\sum_j|d_j|^2-2|\sum_jd_j|^2$.
Thus (V2) removes one factor of the product by differentiation and a second by cancellation. This accounts for the $N-2$ factors that still have to be estimated.

The identity
\[
 \sum_j|1-z_j(h)|^2=N\bigl(|1-s|^2+h^2 V\bigr)
 \le N\bigl(|1-s|^2+1-|s|^2\bigr)=2N(1-\Rea s)
\]
also bounds the sum over any $N-2$ indices. For these indices, write the vector of remaining factors as $(1-t)\mathbf 1+t(1-z_j(h))$. The Euclidean triangle inequality, after division by $\sqrt{N-2}$, bounds its root mean square by $(1-t)+t\sqrt{2\epsilon}=1-Bt$. The arithmetic--geometric mean inequality now gives $|P_{jk}|\le(1-Bt)^{N-2}$, including when one of the factors vanishes. Integration of (V2) with respect to $h$ therefore yields
\[
 \left|F(1)-F(0)\right|\le \tfrac12 NV t^2(1-Bt)^{N-2}.
\]

In detail, for $S=\{1,\ldots,N\}\setminus\{j,k\}$ the two inequalities used above read
\[
 |P_{jk}|^{1/(N-2)}
 \le \left(\frac{1}{N-2}\sum_{\ell\in S}|1-t z_\ell(h)|^2\right)^{1/2}
 \le (1-t)+t\left(\frac{1}{N-2}\sum_{\ell\in S}|1-z_\ell(h)|^2\right)^{1/2}.
\]
The last mean square is at most $2\epsilon$. Equation (V2) therefore bounds $|F'(h)|$ by $hNVt^2(1-Bt)^{N-2}$, and $\int_0^1h\,dh=1/2$ supplies the factor in the estimate for $F(1)-F(0)$.

A change of variable $x=Bt$ and the beta integral give
\[
 \int_0^a t^2(1-Bt)^{N-2}dt\le B^{-3}\int_0^1 x^2(1-x)^{N-2}dx
   = \frac{2}{B^3(N-1)N(N+1)}.
\]

Here $0\le Ba\le1$, so extending the integral in $x$ to $1$ only increases it. Integrating by parts twice gives
\[
 \int_0^1 x^2(1-x)^{N-2}\,dx
 =\frac{2}{N-1}\int_0^1 x(1-x)^{N-1}\,dx
 =\frac{2}{(N-1)N}\int_0^1(1-x)^N\,dx .
\]
The last integral is $1/(N+1)$. Multiplication by $nNV/2$ now leaves precisely $V/[(N-1)B^3]$.

Finally, integrating the common-mean product gives $n\int_0^a(1-ts)^Ndt=[1-(1-as)^n]/s$. Substitution into the preceding bound proves (V1).
\end{proof}
\paragraph{Sharp limiting coefficient.}
For each fixed $N$, the limiting coefficient $1/(N-1)$ in (V1) is sharp as $B\to1$. Put $a=1$ and $y_j=1-\varepsilon+\varepsilon d_j$, where the fixed real $d_j\in[-1,1]$ have mean zero and positive variance. For $0<\varepsilon<1/2$ these points lie in $[0,1]$ and have mean $s=1-\varepsilon$, so $B>0$ for all sufficiently small $\varepsilon$. On expanding the product in powers of the deviations $y_j-s$, the linear term vanishes. The quadratic term is $-NVt^2(1-ts)^{N-2}/2$ because $\sum_{i<j}(y_i-s)(y_j-s)=-NV/2$. Terms of order at least three are $O_N(\varepsilon^3)$ uniformly on $[0,1]$. Integrating the expansion and dividing the signed difference by $V$, we obtain in the limit $\varepsilon\downarrow0$ the value $-nN\int_0^1 t^2(1-t)^{N-2}dt/2=-1/(N-1)$.

For this family one can see the order of the remainder directly by putting $u=1-ts$ and $d=(d_1,\ldots,d_N)$. If $e_k(d)$ is the sum of the products of the $d_j$ over all $k$-element subsets, expansion gives
\[
 \prod_j(u-t\varepsilon d_j)
 = u^N + t^2\varepsilon^2 e_2(d)u^{N-2}
   +\sum_{k=3}^N(-t\varepsilon)^k e_k(d)u^{N-k}.
\]
Here $e_1(d)=0$ and $2e_2(d)=-\sum_jd_j^2$. Since $|t|,|u|,|d_j|\le1$, the last sum has modulus at most $2^N\varepsilon^3$.
As $V=\varepsilon^2 N^{-1}\sum_jd_j^2$, its contribution divided by $V$ tends to zero uniformly. This leaves exactly the integral of the quadratic term computed above.

\begin{corollary}
Let $p$ be monic of degree $n=N+1$, $N\ge3$, with all its zeros in the closed unit disk. Suppose $0<a\le1$ is a marked zero, distinct from every critical point $\zeta_j$, and $N^{-1}\sum_j |a-\zeta_j|^{-2}\le1$. Put $y_j=1/(a-\zeta_j)$ and $q=|1-as|$, with $s$, $V$ and $B$ as in the theorem. If $B>0$, then
\[
 1-a|s| \le q^n + \frac{|s|V}{(N-1)B^3},\qquad 0\le V\le1-|s|^2.
 \tag{V3}
\]
In particular, if $a\ge1/2$ and $2/[(N-1)B^3]\le a/2$, then
\[
 (1-a)+\frac a2(1-|s|)\le q^n,
 \qquad \max_j|y_j-s| \le \sqrt{4N/a}\,q^{n/2}.\tag{V4}
\]
Consequently, along any sequence satisfying the hypotheses with $n\to\infty$, $a\to1$ and $s\to1$, the bounds in (V4) hold for all sufficiently large $n$.
\end{corollary}
\begin{proof}
Factoring $p$ at zero gives $|p(0)|\le a$, while the factorization of $p'$ gives $|p'(a)|=n/\prod_j|y_j|\ge n$ by the arithmetic--geometric mean inequality. Integrating the factorization of $p'$ from zero to $a$, we obtain
\[
 -p(0)=p'(a)\int_0^a\prod_j(1-t y_j)dt.
\]

For the normalization at the marked zero we have, more explicitly,
\[
 \prod_{j=1}^N|y_j|
 \le \left(\frac1N\sum_{j=1}^N|y_j|^2\right)^{N/2}\le1,\qquad
 p'(a-t)=p'(a)\prod_{j=1}^N(1-ty_j).
\]
Writing $I=n\int_0^a\prod_j(1-ty_j)\,dt$, the preceding integral identity yields $|I|=n|p(0)|/|p'(a)|\le |p(0)|\le a$.
If $C=1/[(N-1)B^3]$, (V1) gives the comparison
\[
 1-q^n\le |1-(1-as)^n|\le |s|\,|I|+|s|CV\le a|s|+|s|CV .
\]

It follows that $n$ times the modulus of the integral is at most $a$. Applying (V1), multiplying by $|s|$ and using $|1-(1-as)^n|\ge1-q^n$ proves (V3). The exact identity $V=N^{-1}\sum|y_j|^2-|s|^2$ yields its variance bound. Under the additional hypothesis in (V4), the bound $V\le 2(1-|s|)$ makes the last term in (V3) at most $a(1-|s|)/2$. Subtracting this term gives the first inequality in (V4), and the second follows because each squared deviation is at most $NV$.

In particular, writing $m=|s|$, we have $1-am=(1-a)+a(1-m)$. The first inequality in (V4) therefore gives $1-m\le 2q^n/a$, and hence
\[
 |y_j-s|^2\le\sum_k|y_k-s|^2=NV\le 2N(1-m)\le (4N/a)q^n
 \qquad(1\le j\le N).
\]
For the assertion about sequences, $\epsilon=N(1-\Rea s)/(N-2)\to0$, so that $B\to1$. It follows that $2C\to0$ while $a/2\to1/2$, as required.

\end{proof}

\subsection*{The first coefficient when the mean is close to the boundary}
\begin{corollary}\label{q-small}
Write $s=N^{-1}\sum_j1/(a-\zeta_j)$ and $q=|1-as|$ as above. Under the quadratic reciprocal hypothesis of the preceding corollary, $a=1$ and $p(z)=z^n-1$ whenever $N\ge4$ and
\[
 q\le\min\left\{\frac1{64},R_N\right\},\qquad
 R_N=\frac18\left(1-\frac2N\right)\left(1-\frac4{N+1}\right)^2.
\]
Both factors in the definition of $R_N$ increase with $N$ for $N\ge4$.
In particular $R_4=1/400$ gives the common radius $1/512$ for
$n\ge5$, and $R_7=5/224>1/64$ gives the radius $1/64$ for $n\ge8$.
It follows that no sequence with $\min_j|a-\zeta_j|>1$ can have $n\to\infty$, $a\to1$ and $s\to1$.
\end{corollary}
\begin{proof}

The estimate (V3) controls the loss in modulus $1-a|s|$. To also control the argument of $s$ we use the first coefficient of $p$, which gives a second expression for the mean of its critical points. The two expressions will bound $q$ itself by higher powers of $q$.

If $q=0$, then $a=s=1$. The quadratic hypothesis then gives $V=0$, so every $y_j=1$ and all the critical points are zero. Since $p$ is monic and $p(1)=0$, we obtain $p(z)=z^n-1$. We may therefore assume $q>0$. Put $I=n\int_0^a\prod_j(1-ty_j)dt$, $P=|p(0)|$, $c=a-1/s$, $\tau=q^n$, and write $b_1,\ldots,b_n$ for all source zeros, counting multiplicities. When $I\ne0$, we first show that
\begin{equation*}
 \frac q{|s|}\le \frac{2V}{|s|^2}
        +\frac{-2\log P}{nP}
        +\frac{(q^2+a^2V)^{N/2}}{|I|}. \tag{V5}
\end{equation*}
Indeed $|y_j|\ge1/(1+a)\ge1/2$. Expand the reciprocal about $s$, retaining the exact quadratic remainder. Since the deviations sum to zero, their linear contribution vanishes and we obtain
\[
  \frac1N\sum_j\zeta_j
  =c-\frac1N\sum_j\frac{(y_j-s)^2}{s^2 y_j}.
\]

In fact $1/y_j=1/s-(y_j-s)/s^2+(y_j-s)^2/(s^2y_j)$, and $\zeta_j=a-1/y_j$. The bound $|y_j|^{-1}=|a-\zeta_j|\le a+1$ follows from Gauss--Lucas. Thus, with $\mu=N^{-1}\sum_j\zeta_j$, we have
\[
 |c-\mu|\le\frac1{N|s|^2}\sum_j\frac{|y_j-s|^2}{|y_j|}
 \le \frac{2V}{|s|^2}.
\]
Here $s\ne0$ since $q<1$. To identify $\mu$ directly from $p$, write $p(z)=z^n+A z^{n-1}+\cdots$. The coefficient of $z^{N-1}$ in $p'(z)/n$ is $(N/n)A$, so comparison of coefficients gives $\mu=-A/n=n^{-1}\sum_i b_i$.

Comparison of the first coefficients of $p$ and $p'$ shows that the mean of the zeros equals the mean of the critical points. Moreover, $P=\prod_i|b_i|\le |b_j|$ for every $j$ and $P\ge|I|$ as in (V3). Therefore
\[
 \left|\frac1n\sum_i (1/b_i-\overline{b_i})\right|
 \le \frac{\sum_i(1-|b_i|^2)}{nP}\le\frac{-2\log P}{nP}.
\]

The assumption $I\ne0$ gives $P>0$, so every $b_i$ is nonzero. For $r=|b_i|\le1$ the elementary bound $1-r^2\le-\log(r^2)$ follows from $\log u\le u-1$. Since $1/b_i-\overline{b_i}=(1-|b_i|^2)/b_i$, summing the logarithms gives precisely the factor $-2\log P$. In particular $|\mu|=|\overline\mu|$ is bounded by the displayed error plus $|n^{-1}\sum_i1/b_i|$.

Finally logarithmic differentiation at zero and evaluation of the derivative product give
\[
 \left|\frac1n\sum_i1/b_i\right|
 =\frac{|p'(0)|}{n|p(0)|}
 =\frac{|\prod_j(1-a y_j)|}{|I|}
 \le\frac{(q^2+a^2 V)^{N/2}}{|I|}.
\]

Indeed $p'(0)/p(0)=-\sum_i1/b_i$, and $p'(0)=p'(a)\prod_j(1-a y_j)$, while $nP=|p'(a)|\,|I|$. For the resulting product the relevant second moment is
\[
 \frac1N\sum_j|1-a y_j|^2
 =\frac1N\sum_j|(1-as)-a(y_j-s)|^2
 =|1-as|^2+a^2V,
\]
the mixed term having zero mean. Applying the arithmetic--geometric mean inequality to these squared moduli explains the exponent $N/2$.

The last estimate follows by taking the arithmetic mean of the squared moduli and using the fact that the deviations sum to zero. Combining these statements with $|c|=q/|s|$ proves (V5).

Write $C=1/[(N-1)B^3]$. We have $a,|s|,\Rea s\ge1-q$ because $\Rea(as)\ge1-q$.
Let $b=(N+5)/(2(N+1))$. The bound $q\le R_N$ implies
$\sqrt{2N(1-\Rea s)/(N-2)}\le(N-3)/(2(N+1))=1-b$.
Thus $B\ge b>0$. To see that $(N-1)b^3\ge2187/1000$, clear denominators and set $t=N-4$. The resulting difference
$1000(N-1)(N+5)^3-2187[2(N+1)]^3$ equals
$t(1000t^3+12504t^2+61560t+145800)$.
Consequently $C\le1000/2187<1/2$ and
$a-2C\ge63/64-2000/2187>1/16$.

To see how the threshold was chosen, substitute $1-\Rea s\le q\le R_N$ to obtain
\[
 \frac{2N}{N-2}(1-\Rea s)\le \frac{2N}{N-2}R_N
 =\frac14\left(1-\frac4{N+1}\right)^2.
\]
The factor in parentheses is positive for $N\ge4$, so taking square roots gives the stated lower bound on $B$.

For the remainder put $m=|s|$. Since $mV\le2(1-m)$, (V3) implies
\[
 (1-a)+(a-2C)(1-m)\le\tau .
\]
In particular $1-m<16\tau$ for $\tau>0$, and $V\le2(1-m)<32\tau$. In (V1) the mean product has modulus at least $(1-\tau)/m\ge1-\tau$, and the error satisfies $CV<16\tau$. This bounds $P$ and $|I|$ away from zero before either is used as a denominator.

Equation (V3) now gives $(1-a)+(1-|s|)/16\le\tau$ and $V<32\tau$. It follows from (V1) that $P\ge|I|\ge(1-\tau)/|s|-CV>1-17\tau>1/2$, so we may apply (V5). Its first two terms are at most $67\tau$ and $28\tau$ respectively, since $|s|\ge63/64$, $-\log P\le 2(1-P)$ for $P>1/2$, and $n\ge5$. Since $q^2+32q^n<4q^2$, (V5) gives
\[
 q\le 95q^n+2(2q)^N.
\]

For the constants in this last estimate, $2V/m^2<64(64/63)^2\tau<67\tau$. Also $-\log P=\int_P^1 dx/x\le2(1-P)$, so that
\[
 \frac{-2\log P}{nP}\le \frac{8(1-P)}n < \frac{136}{5}\tau<28\tau .
\]
Finally $q^2+32q^n=q^2(1+32q^{n-2})<4q^2$ because $n\ge5$ and $0<q\le1/64$. The last term in (V5) is therefore at most $2(2q)^N$, and its left-hand side $q/m$ is at least $q$. After division by $q$, the two terms on the right are $95q^{n-1}$ and $4(2q)^{N-1}$, respectively.

Dividing by $q$ gives a contradiction, since $95/64^4+4/32^3<1$. The assertion about sequences follows by observing that $q\to0$ whenever $a\to1$ and $s\to1$.
\end{proof}

\section{A complex integration path and a trace bound}\label{sec:contour}
The real interval suffices once $as$ is close to one. To determine when this happens, we deform the interval into two segments meeting at the center of the disk where the product estimate gives the greatest decay. The normalized trace determines this center. We first estimate $|s|$ and then use a coefficient of $p$ to control its phase.
The same estimate must serve two purposes. When $N\alpha$ is large, we need a small numerical constant multiplying $V$. Later we shall have only a lower bound on $N|s|^2$, and we will need to know how this constant changes with $N|s|^2$. We keep that dependence in the following integral.

\begin{lemma}[Uniform phase bound]\label{phase-constants}
For $\phi(u)=u-u^2/2$ and
$B(L)=L^2\int_0^1u^2 e^{-L\phi(u)}du$, the following bound holds
\[
 L B(L)\le 4+\frac{18}{L}\qquad(L\ge9).
\]
\end{lemma}
\begin{proof}
Set $x=L/2$, $T_j(x)=x^j/j!$ and
$M_j=2/[(2j+1)(2j+2)(2j+3)]$. The identity $1/2-\phi(u)=(1-u)^2/2$ allows us to expand the integrand in powers of $(1-u)^2$. Since $\int_0^1u^2(1-u)^{2j}du=M_j$, the exponential series gives
\[
 x e^x L B(L)=8x^4\sum_{j\ge0}M_j T_j(x).
\]
Indeed, setting $w=1-u$ gives
\[
 M_j=\int_0^1(1-w)^2w^{2j}\,dw
     =\frac1{2j+1}-\frac2{2j+2}+\frac1{2j+3}.
\]
The expansion in the display above now follows from
$e^{-L\phi(u)}=e^{-x}\exp(x(1-u)^2)$ and $xL^3=8x^4$.
Termwise integration is justified by monotone convergence, since all terms are nonnegative. To estimate the tail, observe that for $j\ge6$
\[
 8(j+2)(j+3)(j+4)\le [4(j+4)+9](2j+1)(2j+3).
\]
Indeed, on writing $j=6+z$ the right-hand side exceeds the left by
$8z^3+204z^2+1588z+3795$. Cancelling factorials gives
$8x^4M_jT_j\le4xT_{j+3}+9T_{j+4}$.

To see the comparison at a single index, divide both sides by $T_{j+4}(x)$. On the left this gives $8M_j(j+1)(j+2)(j+3)(j+4)$, while the right side becomes $4(j+4)+9$. Substituting $M_j$ cancels $2(j+1)$ against $2j+2$ and gives precisely the preceding polynomial inequality. Summing from $j=6$ consequently uses the terms beginning with $T_9$ in $4x e^x$ and with $T_{10}$ in $9e^x$.

For the first six terms define
\[
 P(x)=4x\sum_{j=0}^8T_j(x)+9\sum_{j=0}^9T_j(x)
       -8x^4\sum_{j=0}^5M_jT_j(x).
\]
To check the sign of this polynomial of degree nine, write it in the form
\[
 \begin{aligned}
 P(x)={}&9+13x+\tfrac{17}2x^2+\tfrac72x^3-\tfrac{13}8x^4-\tfrac1{40}x^5\\
       &+\tfrac{13}{1680}x^6+\tfrac{31}{15120}x^7
         +\tfrac{457}{1330560}x^8+\tfrac{89}{1921920}x^9 .
 \end{aligned}
\]
For an explicit check of the shift, write $P(x)=\sum p_jx^j$ and $x=\tfrac92(1+w)$. The coefficient of $w^k$ is $\sum_{j=k}^9\binom jk(9/2)^jp_j$. For $k\ge6$ every summand is positive. In increasing order of $k$ the first six coefficients are bounded below by $79,162,2120,7336,11206,9597$, respectively, on substituting the fractions just displayed. Since $w=(2x-9)/9$, this calculation also checks the signs for the shift $y=2x-9$.

Substitution of $x=(9+y)/2$ gives positive coefficients in every degree, so $P(x)\ge0$ for $x\ge9/2$. Combining this bound for the first six terms with the termwise estimate for the tail yields $x e^x L B(L)\le(4x+9)e^x$, as required.
\end{proof}

\begin{theorem}
Let $p$ be a monic polynomial of degree $n\ge2$ with all its zeros in the closed unit disk, and let $0<a\le1$ be a simple zero. Denote the $N=n-1$ critical points by $\zeta_j$, counted with multiplicity, and put $y_j=1/(a-\zeta_j)$ and $s=N^{-1}\sum_j y_j$. If
\[
     \frac1N\sum_j |y_j|^2\le1,\qquad
     N\alpha\ge20,\qquad \alpha=2a\Rea s-a^2,
\]
then $a=1$ and $p(z)=z^n-1$. In particular, any counterexample to Sendov's bound at the marked zero would satisfy
\[
      (n-1)\bigl(2a\Rea s-a^2\bigr)<20.
\]
Thus for any fixed $\eta>0$, the lower bound $\alpha\ge\eta$ rules out such a counterexample in every degree $n\ge1+20/\eta$.
\end{theorem}
\begin{proof}
Set $m=|s|$, $x=as$, $q=|1-x|$, $V=N^{-1}\sum|y_j-s|^2$, $v=N\alpha$. The quadratic assumption gives $V\le1-m^2$. Since $\alpha=m^2-|a-s|^2\le1$ and $v\ge20$, one has $N\ge20$, $m>0$ and $0<\alpha\le m^2$. Also
\[
 q^2+a^2 V =N^{-1}\sum_j|1-a y_j|^2\le1-\alpha.
 \tag{U1}
\]
In fact $1-a y_j=(1-as)-a(y_j-s)$, and averaging its squared modulus removes the cross term. Expanding the same average about zero instead gives $1-2a\Rea s+a^2N^{-1}\sum_j|y_j|^2$, which is at most $1-2a\Rea s+a^2$.
We write the integral from the marked root to the origin in scaled coordinates and compare it with the integral obtained by replacing every $y_j$ by $s$. Denote these two quantities by
\[
 J=n\int_0^1\prod_j(1-t a y_j)dt,\qquad J_0=\frac{1-(1-x)^n}{x}.
\]

Here the substitution $z=a(1-t)$ gives $p'(z)/p'(a)=\prod_j(1-ta y_j)$, whereas integration of $(1-tx)^N$ gives $n\int_0^1(1-tx)^Ndt=[1-(1-x)^n]/x$. Thus in these coordinates $xJ_0$ differs from $1$ by a term of modulus $q^n$. It will suffice to estimate the error made in passing from $xJ_0$ to $xJ$.

Since $-p(0)=(a/n)p'(a)J$, $|p(0)|\le a$, and $|p'(a)|/n=1/\prod|y_j|\ge1$, we have $|J|\le1$. Both integrands are polynomials in $t$, so we may integrate along the two segments $0\to C\to1$, where
\[
 C=\overline{s}/a,\qquad \gamma(t)=2\Rea(tx)-a^2|t|^2=m^2-a^2|t-C|^2 .
\]
In particular $\gamma(1)=\alpha$, and $\gamma$ attains its maximum at $C$.

\paragraph{From zero to $C$.} Put $z_j(h)=s+h(y_j-s)$ and $F(h,t)=\prod_j(1-ta z_j(h))$. Equation (V2) of Section~2 gives
\[
 \partial_h F=-\frac{h a^2t^2}{2N}\sum_{j\ne k}(y_j-y_k)^2
       \prod_{\ell\ne j,k}(1-ta z_\ell(h)).
\]
The identity is algebraic and therefore still holds for complex $t$. Since the deviations from $s$ sum to zero, $N^{-1}\sum|z_j(h)|^2=m^2+h^2 V\le1$. Hence $\sum_j|1-ta z_j(h)|^2\le N(1-\gamma(t))$. Dropping two terms preserves this upper bound for the sum. Applying the arithmetic--geometric mean inequality to the remaining $N-2$ squared moduli and then using $1+u\le e^u$, we obtain
\[
 \left|\prod_{\ell\ne j,k}(1-ta z_\ell(h))\right|
 \le\left(\frac{N(1-\gamma(t))}{N-2}\right)^{(N-2)/2}
 \le \exp(1-N\gamma(t)/2).
\]

The same bound holds by continuity if the product vanishes. Write $\lambda=Nm^2\ge v$. On $t=uC$, $0\le u\le1$, one has $\gamma(t)=m^2(2u-u^2)$. Using $\sum|y_j-y_k|^2=2N^2V$ and $|x|a^2|C|^3=m^4$, we integrate first in $h$ and then in $u$ to bound the contribution of the first segment to $|x(J-J_0)|$ by
\[
 \frac{en}{2N} V\lambda^2\int_0^1u^2e^{-\lambda(u-u^2/2)}du.
\]

To derive the integrand, first write the base in the bound for the product of the remaining $N-2$ factors as $1+(2-N\gamma(t))/(N-2)$. Raising the exponential bound to the power $(N-2)/2$ then produces the exponent $1-N\gamma(t)/2$. Integrating the identity for $\partial_h F$ from $h=0$ to $1$ yields
\[
 |F(1,t)-F(0,t)|
 \le \frac{a^2|t|^2}{2N}\, e^{1-N\gamma(t)/2}
      \left(\sum_{j\ne k}|y_j-y_k|^2\right)\int_0^1 h\,dh.
\]
The sum over ordered pairs is $2N\sum_j|y_j|^2-2|\sum_jy_j|^2=2N^2V$. Hence the right side is $eNa^2V|t|^2 e^{-N\gamma(t)/2}/2$. In particular the difference of the two products vanishes to second order at the initial endpoint, where exponential decay alone gives no gain.

In passing from this pointwise bound to the first contour integral, the powers of $|C|$ come from $|t|^2|dt|=|C|^3u^2du$ and the factor $n|x|$ in front of the integral. Their combined coefficient is $enN V m^4/2=enV\lambda^2/(2N)$. We shall use this expression again with a smaller lower bound on $\lambda$ in Section~\ref{sec:large}.

By Lemma~\ref{phase-constants}, $\lambda B(\lambda)\le4+18/\lambda
\le49/10$ for $\lambda\ge20$. With $e<11/4$ and $n/N\le21/20$
the first segment therefore contributes at most
\[
  \frac{en}{2N}B(\lambda)V
  \le \frac{11}{4}\frac{21}{20}\frac{49}{10}\frac1{40}\,V
  <0.39V\quad\text{if } V>0. \tag{U2}
\]
When $V=0$ the contribution vanishes by (V2).

\paragraph{From $C$ to one.} Keeping all $N$ factors in the arithmetic--geometric mean inequality gives $|F(1,t)|,|F(0,t)|\le e^{-N\gamma(t)/2}$. Write $d=m^2-\alpha=a^2|1-C|^2$. On $t=C+u(1-C)$ one has $\gamma(t)=\alpha+d(1-u^2)\ge\alpha+d(1-u)$. If $d=0$ this segment has zero length. Otherwise its contribution to $|x(J-J_0)|$ is at most
\[
  2n m\sqrt d\,e^{-v/2}\min\bigl(1,2/(Nd)\bigr)
 \le 2(n/N)\sqrt{2(v+2)}\,e^{-v/2}.\tag{U3}
\]
Indeed, for $h=0,1$ the bound with all factors follows from $(1-\gamma(t))^{N/2}\le e^{-N\gamma(t)/2}$. The triangle inequality contributes a factor of two along this segment, and multiplying its length by $|x|$ gives $am\,|1-C|=m\sqrt d$. The remaining integral is
\[
 \int_0^1 e^{-Nd(1-u)/2}\,du
  =\frac{1-e^{-Nd/2}}{Nd/2}
  \le \min\{1,2/(Nd)\}.
\]
For small $Nd$ the first bound retains the length of the segment. As $Nd$ grows the second bound also takes account of decay away from its endpoint.

To verify the last inequality, put $L=Nd>0$. The coefficient of $e^{-v/2}$ on the left is then $2(n/N)\sqrt{(v+L)L}\min(1,2/L)$. This expression increases up to $L=2$ and decreases thereafter, which proves (U3).

Indeed, after squaring and removing the positive constant \(4(n/N)^2\), the expression is \((v+L)L\) for \(0<L\le2\) and \(4(1+v/L)\) for \(L\ge2\). The two values agree at \(L=2\), where they equal \(2(v+2)\).

\paragraph{The distance from $as$ to one.} Put $D=1-|x|$ and $\delta=(1+2(n/N)\sqrt{2(v+2)})e^{-v/2}$. Since $q^n\le e^{-v/2}$ by (U1) and $|J|\le1$, (U2)--(U3) imply
\[
 D\le 0.39 V+\delta,\qquad
   |x|\,|J|\ge 1-0.39V-\delta.
\]
More explicitly, (U1) gives $q^n\le(1-\alpha)^{n/2}\le e^{-n\alpha/2}\le e^{-v/2}$. Subtracting $xJ_0$ inside $1-xJ$ now gives
\[
 1-|x|\le 1-|xJ|\le |1-xJ|
 \le q^n+|x(J_0-J)|\le 0.39V+\delta .
\]
The first inequality uses $|J|\le1$. To absorb the variance, use $|x|=am\le m$ and $V\le(1-m)(1+m)\le2D$. It follows that $(1-2\cdot0.39)D\le\delta$, with $1-2\cdot0.39=11/50$. Substituting this bound back into $1-|xJ|\le0.39V+\delta$ also gives $1-|J|\le 50\delta/11$.

Thus $V\le2(1-m)\le2D$ gives $D\le50\delta/11$, $V\le100\delta/11$, and $|J|\ge1-50\delta/11$. To bound $\delta$, observe that $n/N\le21/20$ and that $\sqrt{2(v+2)}e^{-v/2}$ decreases for $v\ge20$. Evaluating at $v=20$ gives $\delta<1/1450$, since $\sqrt{44}<20/3$ and $e^{10}>21750$. The latter inequality follows from $e>2717/1000$ by retaining the first seven terms of the exponential series at one. In particular $a,m,|J|,a|J|,P>0.99$, with $P=|p(0)|\ge a|J|\ge1-100\delta/11$.

For the product bound in the last assertion, each of $a$ and $m$ is at least $am=1-D$. Hence $a|J|\ge(1-50\delta/11)^2\ge1-100\delta/11$, and $P=(|p'(a)|/n)a|J|\ge a|J|$. These estimates control the moduli. The phase still enters through $q=|1-as|$, since a lower bound on $|as|$ alone gives no upper bound on this distance.

It remains to control the phase of $s$. In the present coordinates, the estimate (V5) obtained from the first coefficient reads
\[
 \frac q m\le \frac{2V}{m^2}+\frac{-2\log P}{nP}
        +\frac{(q^2+a^2V)^{N/2}}{a|J|}. \tag{U4}
\]
The denominator in the last term is the modulus of $I=aJ$, obtained by scaling the real interval in (V5). For the numerical comparison, $V\le100\delta/11$ and $1-P\le100\delta/11$ give
\[
 \frac{2V}{m^2}
  \le \frac{200\delta}{11(0.99)^2},\qquad
 \frac{-2\log P}{nP}
  \le \frac{2(1-P)}{nP^2}
  \le \frac{200\delta}{11\cdot21(0.99)^2}.
\]
In the middle inequality we used $-\log P=\int_P^1 du/u\le(1-P)/P$ for $0<P\le1$. Finally $e^{-v/2}\le\delta$ from its definition.

Here we use the equality of the means of the zeros and critical points, logarithmic differentiation at zero, and the bound $|y_j|\ge1/2$ from the Gauss--Lucas theorem. Equation (U1) bounds the numerator of the last term by $e^{-v/2}$. Because $m,P,a|J|>0.99$, the first term is less than $19\delta$. The second is at most $2(1-P)/(nP^2)<\delta$, using $n\ge21$. The third is less than $2\delta$. Thus $q<22\delta<1/64$. Since $n\ge21$, the radius $1/64$ in Corollary~\ref{q-small} applies and gives $a=1$ and $p(z)=z^n-1$.
\end{proof}

\section{Degrees 257 and above}\label{sec:large}
Reflection gives a lower bound on $v$ that increases with the degree, while the contour argument gives the upper bound $v<20$. We use the phase estimate once more to exclude the values that could satisfy both conditions.

\begin{theorem}\label{large}
For $n\ge257$ the quadratic reciprocal hypothesis~\eqref{hyp}
cannot hold at an interior marked zero $0<a<1$.
\end{theorem}
\begin{proof}
Suppose $N=n-1\ge256$. Equation~\eqref{polar-exp} implies
$v>2\log(N/2)\ge2\log 128>9.7$. We also have $a<4/5$, since otherwise
$0<Y\le9N/50$ and convexity would imply
$1\le F_\alpha(9N/50)<e^{9v/50}/(9N/50)<1$,
since $v<20$, $e^{18/5}<40$ and $9N/50>40$.

Here the two uses of reflection concern different values of its argument. At $N/2$, inequality~\eqref{polar-exp} yields $e^{v/2}\ge (N/2)(1+\alpha)+e^{-N/2}>N/2$. At the other value $Y_*=9N/50$, the assumption $a\ge4/5$ gives $Y=N(1-a^2)/2\le Y_*$. Convexity with $F_\alpha(0)=1$ shows
\[
 1\le F_\alpha(Y)\le (1-Y/Y_*)+(Y/Y_*)F_\alpha(Y_*),
\]
so $F_\alpha(Y_*)\ge1$. On the other hand direct integration gives $F_\alpha(Y_*)=(e^{\alpha Y_*}-e^{-Y_*})/[Y_*(1+\alpha)]$. The contour theorem excludes $v\ge20$ at an interior zero. Consequently $\alpha Y_*=9v/50<18/5$, whereas $Y_*\ge9\cdot256/50>40$. This proves the claimed restriction on $a$.

For completeness, the elementary exponential bounds needed at these endpoints can be checked with the same finite sum. Let $S_{10}(u)=\sum_{j=0}^{10}u^j/j!$. For $0<u<12$, the ratio between successive terms after the first omitted term is at most $u/12$, so
\[
 S_{10}(u)<e^u\le S_{10}(u)+\frac{u^{11}}{11!\,(1-u/12)}.
\]
At $u=97/20$ the lower bound exceeds $125$ and the upper bound is less than $128$. At $u=18/5$ the upper bound is less than $37$.

Put $L=Nm^2\ge v$ and $A(L)=e(N+1)B(L)/(2N)$. Combining the estimates on the two segments of the contour in Section~\ref{sec:contour} gives
\[
 1-am\le A(L)V+\delta,\qquad
 \delta=\left(1+2(n/N)\sqrt{2(v+2)}\right)e^{-v/2}<9/100 .
\]

The coefficient on the first segment here is the expression preceding (U2), before the specialization to $L\ge20$. That expression and (U3) only require $N>2$, $m>0$ and $\alpha\ge0$. Substituting both in $1-am\le q^n+|x(J-J_0)|$ gives the displayed inequality throughout the present range $97/10<v<20$.

For the bound on $\delta$, both terms decrease for $v\ge9.7$, and at $v=9.7$ we may use $e^{9.7/2}>125$, $\sqrt{2(9.7+2)}<5$ and $n/N\le257/256$.

Indeed the logarithmic derivative of $\sqrt{2(v+2)}e^{-v/2}$ is $1/[2(v+2)]-1/2<0$. It follows that $\delta<(1+10\cdot257/256)/125<9/100$. Meanwhile $L\ge v>9$ gives $L B(L)\le4+18/L\le6$, so the factor $e(N+1)/(2N)$ in $A(L)$ converts the phase bound into an upper bound on $AL$.

Set $A=A(L)$. By Lemma~\ref{phase-constants},
$AL\le41/5$, because $e<68/25$ and
$e(N+1)/(2N)\le41/30$. If $A\le2/3$, then
\[
 1-am-AV\ge 1-\tfrac45m-\tfrac23(1-m^2)
 =\tfrac7{75}+\tfrac23(m-\tfrac35)^2>\tfrac9{100}.
\]
This contradicts the upper bound $1-am-AV\le\delta$. It follows that $A>2/3$ and $L<123/10<13$, whence $m^2<13/256$.
As $L>97/10$, we also have $A\le17/20$.

The two bounds on $L$ now have complementary uses. From $A>2/3$ one gets $L\le(41/5)/A<123/10$, so the mean square $m^2=L/N$ is small. Conversely $L>97/10$ bounds $A$ itself by $(41/5)/(97/10)=82/97<17/20$. The estimate that follows uses both restrictions in $am+AV$, together with $V\le1-m^2$.

The relation $\alpha=2a\Rea s-a^2>0$ implies $a\le2m$, giving
\[
 1-am-AV\ge\tfrac3{20}-\tfrac{23}{20}m^2
 >\tfrac{469}{5120}>\tfrac9{100}.
\]
This again contradicts $1-am-AV\le\delta$, exhausting the possible values of $L$.
\end{proof}

\section{Polynomial moment bounds for degrees five through 256}\label{sec:finite}
For these degrees we retain the powers in the two contour integrals. Together with reflection they give necessary inequalities in the real parameters $a$ and $m$. We show that on each rectangle of a subdivision of their domain at least one of these inequalities fails throughout the corresponding interval of integer degrees.

Let $p$ be monic with all zeros in the closed unit disk, and rotate a
simple interior zero to $0<a<1$. Write $c_1,\ldots,c_N$ for the
critical points counted with multiplicity, where $N=n-1$, and let $\E$ denote the average over $j=1,\ldots,N$. Under the
assumption $N^{-1}\sum|a-c_j|^{-2}\le1$, we have
\[
 \begin{gathered}
 y_j=(a-c_j)^{-1},\quad s=\E y,\quad m=|s|,\quad V=\E|y-s|^2\le1-m^2,\\
 \alpha=2a\Re s-a^2,\qquad b=1-a^2,\qquad q=|1-as|.
 \end{gathered}
\]
Averaging $2a\Re y_j-a^2=b(1-|y_j|^2)+(1-|c_j|^2)|y_j|^2$ gives
$\alpha\ge0$, while $\Re s\le m$ gives
$\alpha\le m^2-|a-m|^2\le1$. In particular $m>0$, and
$q^2=1-\alpha-a^2(1-m^2)$.
\begin{theorem}
Let $p$ be a monic polynomial of degree $5\le n\le256$ with all zeros in the closed unit disk. At each simple zero $a$ with $0<|a|<1$, its critical points satisfy $N^{-1}\sum|a-c_j|^{-2}>1$, where $N=n-1$.
\end{theorem}

We shall prove the theorem by subdividing the range of degrees and
the unit square for $(a,m)$ into boxes with rational endpoints. We
first derive the inequality that will exclude most of these boxes.

\paragraph{Integration from the marked zero to the origin.}
Let $J=n\int_0^1\prod(1-at y_j)dt$, and let $J_0$ be the same
integral with each $y_j$ replaced by $s$. The product formulas for
$p$ and $p'$ give $|J|\le1$, while direct integration gives
$asJ_0=1-(1-as)^n$. To bound the difference, we integrate along the
segments $0\to C\to1$ with $C=\bar s/a$. Writing
$C_N=(N/(N-2))^{(N-2)/2}$, we obtain the necessary inequality
\begin{align}
 1-am&\le q^n+E_0+E_1,\label{finiteineq}\\
 E_0&=\frac{nN}2 C_N\,m^4(1-m^2)
  \int_0^1 u^2[1-m^2(2u-u^2)]^{(N-2)/2}du,\notag\\
 E_1&=2n m\sqrt{m^2-\alpha}
   \int_0^1[1-m^2+(m^2-\alpha)u^2]^{N/2}du.\notag
\end{align}

Indeed, for $0\le h\le1$ the interpolation $z_j(h)=s+h(y_j-s)$
satisfies $\E|z(h)|^2=m^2+h^2V\le1$. For
$F(h,t)=\prod(1-at z_j(h))$, subtracting the products obtained by
omitting the $j$th and the $k$th factors gives
\[
 \partial_h F=-\frac{h a^2t^2}{2N}\sum_{j\ne k}
 (y_j-y_k)^2\prod_{l\ne j,k}(1-at z_l(h)).
\]
The corresponding sum of the squared moduli $|y_j-y_k|^2$ equals
$2N^2V$. On the first segment, where $t=uC$, the mean of the squared
moduli of all $N$ factors is at most $1-m^2(2u-u^2)$. Omitting
two terms can only decrease their sum. Applying the arithmetic--geometric
mean inequality to the remaining $N-2$ squared moduli gives the factor
$C_N$ in the bound for the product. Integration in $h$ and $u$, together
with $|as|a^2|C|^3=m^4$, now gives $E_0$.

More explicitly, the same square completion works on both segments.
Using $\E z(h)=s$ and $\E|z(h)|^2\le1$ gives
\[
 \E|1-at z(h)|^2\le 1-2\Rea(at s)+|at|^2
                      =1-m^2+|at-\bar s|^2 .
\]
On the first segment $at=u\bar s$, so the right side is
$W_0(u)=1-m^2(2u-u^2)$. The arithmetic mean of the remaining $N-2$ squared moduli
is at most $NW_0(u)/(N-2)$. Thus the pairwise identity above gives
\[
 |\partial_hF(h,uC)|\le hNa^2u^2|C|^2 V C_N W_0(u)^{(N-2)/2}.
\]
Integrating $h$ from zero to one supplies the factor $1/2$. The
remaining factor from the contour integral is $n|as|\,|C|$, and
$V\le1-m^2$.

On the second segment, $t=C+u(1-C)$, the mean squared modulus is at
most $1-m^2+(m^2-\alpha)u^2$ for the factors in each of $F(1,t)$
and $F(0,t)$. Adding the two product bounds and using
$a|1-C|=\sqrt{m^2-\alpha}$ gives $E_1$. Thus
$|as(J-J_0)|\le E_0+E_1$. Substituting the formula for $asJ_0$
and using $|asJ|\le am$ gives~\eqref{finiteineq}. The argument
uses the mean squared modulus throughout and allows individual
$|y_j|$ to exceed one.

Indeed, on this segment $at-\bar s=u(a-\bar s)$ and
$|a-\bar s|^2=m^2-\alpha$. Its length in the $t$-plane is $|1-C|$,
so the two product bounds enter with the common total factor
$2n|as|\,|1-C|=2nm\sqrt{m^2-\alpha}$. The final triangle inequality
is applied to the identity $1=(1-as)^n+asJ+as(J_0-J)$.

\begin{lemma}[Polynomial moments]\label{lem:moments}
For nonnegative $A,B,C$, put
\[
 Q(t)=A(1-t)^2+B t(1-t)+C t^2,\qquad
 M_{j,k}(Q)=\int_0^1t^jQ(t)^k\,dt\quad (j,k\in\mathbb Z_{\ge0}).
\]
If $(A+Bz+Cz^2)^k=\sum_{h=0}^{2k}d_h z^h$, then
\begin{equation}\label{eq:moment-sum}
 M_{j,k}(Q)=\sum_{h=0}^{2k}
 \frac{d_h}{(2k+j+1)\binom{2k+j}{h+j}} .
\end{equation}
If $0\le Q\le1$ on $[0,1]$ and $p\ge0$ is real, then
$\int_0^1 t^jQ^p\le M_{j,\lfloor p\rfloor}(Q)$.
For any $Q\ge0$ and any integer $k\ge p>0$, the inequality
$M_{0,k}(Q)<1$ implies $\int_0^1Q^p<1$.
\end{lemma}
\begin{proof}
In the coefficient expansion, substitute $z=t/(1-t)$ and multiply by
$(1-t)^{2k}$. The resulting polynomial identity extends to $t=1$.

In particular the coefficient of $z^h$ contributes precisely
$d_h t^h(1-t)^{2k-h}$ to $Q(t)^k$.
Integrating each term against $t^j$ gives~\eqref{eq:moment-sum} since
$\int_0^1t^{j+h}(1-t)^{2k-h}dt=(j+h)!(2k-h)!/(2k+j+1)!$.
When $0\le Q\le1$, the second assertion holds pointwise, with
$Q^0=1$. For the last assertion, Jensen's inequality for the convex
function $u\mapsto u^{k/p}$ gives
$(\int_0^1 Q^p)^{k/p}\le\int_0^1 Q^k$.
\end{proof}

For the weight in~\eqref{eq:moment-sum}, let \(I_{r,s}=\int_0^1t^r(1-t)^s\,dt\). Integration by parts gives \(I_{r,s}=r I_{r-1,s+1}/(s+1)\) for \(r\ge1\), and \(I_{0,s}=1/(s+1)\). Taking \(r=h+j\) and \(s=2k-h\), these identities give
\[
 I_{h+j,2k-h}
 =\frac{(h+j)!\,(2k-h)!}{(2k+j+1)!}
 =\frac{1}{(2k+j+1)\binom{2k+j}{h+j}}.
\]
This explains why each coefficient can be integrated using an integer denominator.

The same moments can be obtained recursively by setting
$T_{k,i,j}=\int_0^1t^i(1-t)^jQ(t)^kdt$ and expanding one factor of
$Q$. This gives
\[
 T_{k+1,i,j}=A T_{k,i,j+2}+B T_{k,i+1,j+1}+C T_{k,i+2,j},
 \qquad T_{0,i,j}=\frac{1}{(i+j+1)\binom{i+j}{i}} .
\]
The recurrence uses only multiplication and addition after the initial
values have been set. In Lean its agreement with the integral is
expressed by \texttt{rationalMoment\_correct}. For the calculation below
we use~\eqref{eq:moment-sum}, which involves $2k+1$ coefficients.

\paragraph{Reflection on a box.}
Consider the intervals $a_l\le a\le a_h$, $m_l\le m\le m_h$ and
$n_l\le n\le n_h$, and choose $0\le b_l\le b\le b_h\le1$ for
$b=1-a^2$. If $\alpha\le x$ for some $0\le x\le1$, then the
quadratic in~\eqref{polar} satisfies
\[
 1+b((1+\alpha)t-1)-b^2t(1-t)
 \ \le\ 1+bH(t),\qquad
 H(t)=(1+x)t-1-b_l t(1-t).
\]

The difference between the right and left sides is
$b(x-\alpha)t+b(b-b_l)t(1-t)$, a sum of two nonnegative terms.
Also $1+H(t)=t(1+x-b_l+b_l t)\ge0$.
Since $H(t)\ge-1$ on the unit interval, $1+zH(t)\ge0$ for
$0\le z\le b_h$. The function
$z\mapsto\int_0^1(1+zH)^{N/2}$ is convex, takes the value one
at zero and, by reflection, is at least one at $z=b>0$.
Convexity therefore gives the same lower bound at $b_h$. In the basis
of Lemma~\ref{lem:moments}, the coefficients of $Q=1+b_hH$ are
\[
 A=1-b_h,\quad B=2(1-b_h)+b_h(1+x-b_l),\quad C=1+b_hx .
\]

They are all nonnegative. For $k=\lfloor n_h/2\rfloor
=\lceil(n_h-1)/2\rceil$, equation~\eqref{eq:moment-sum} gives
one rational comparison for all degrees in the box. If
$M_{0,k}(Q)<1$, Lemma~\ref{lem:moments} contradicts the lower
bound at $b_h$, so $\alpha>x$. We first test $x=1$, since success
would exclude the entire box. Otherwise we bisect $[0,1]$ twenty
times, raising the left endpoint only when the comparison gives
$\alpha>x$. The final left endpoint is a lower bound
$\underline\alpha\le\alpha$.

To see the endpoint argument explicitly, put
$\Phi(z)=\int_0^1(1+zH(t))^{N/2}dt$ and $\theta=b/b_h$. Then
$0<\theta\le1$, and
\[
 1\le\Phi(b)\le (1-\theta)\Phi(0)+\theta\Phi(b_h)
                 =1-\theta+\theta\Phi(b_h).
\]
Thus $\Phi(b_h)\ge1$ even though $H$ may change sign on the interval.

For a concrete box from the subdivision, take $5\le n\le8$ and
$0\le a,m\le1/2$. The geometric bound gives $\alpha\le m^2\le1/4$.
In the reflection test choose $x=1/2$ and use $b_l=3/4$, $b_h=1$.
Then $k=4$ and $Q(t)=3t(1+t)/4$, so direct expansion gives
\[
 M_{0,4}(Q)=\left(\frac34\right)^4
       \sum_{j=0}^4\binom4j\frac1{5+j}
       =\frac{13239}{17920}<1 .
\]
Although $Q(1)=3/2$, its fourth moment is strictly below one.
Reflection would therefore require $\alpha>1/2$, whereas the geometric
bound gives $\alpha\le1/4$ throughout the box. This one calculation
excludes all four integer degrees on the entire rectangle.

\paragraph{The two integrals on the path to the origin.}
For a box that remains after the reflection test, put $z=m_l^2$
and $w=1-z$. The two quadratics in~\eqref{finiteineq} are bounded
above by
\[
 Q_0=1-z(2u-u^2),\qquad
 Q_1=1-z(1-u^2)-\underline\alpha u^2 ,
\]
respectively. Their coefficient triples in the basis of
Lemma~\ref{lem:moments} are $(1,2w,w)$ and
$(w,2w,1-\underline\alpha)$. The first quadratic is affine in
$(1-u)^2$ and the second in $u^2$. Since their endpoint values lie
in $[0,1]$, so do both quadratics on the whole interval. The integrals
in $E_0$ and $E_1$ therefore have the respective upper bounds
\[
 M_{2,\lfloor(n_l-3)/2\rfloor}(Q_0),
 \qquad M_{0,\lfloor(n_l-1)/2\rfloor}(Q_1).
\]
Together with the reflected integral, this gives three applications
of the same moment formula. A sharper estimate is available when the
exponent at the lower degree is a half-integer. For integers
$k,j\ge0$, Cauchy--Schwarz gives
$\int_0^1u^jQ^{k+1/2}\le\sqrt{M_{j,k}(Q)M_{j,k+1}(Q)}$,
which also holds when either moment vanishes. The integer powers displayed above suffice for the subdivision used in
the proof.

One can read the two bounds directly from the differences of the
quadratics. In the first integral the difference between $Q_0$ and
the original base is $(m^2-z)(2u-u^2)$. In the second it is
\[
 Q_1-\bigl[1-m^2+(m^2-\alpha)u^2\bigr]
   =(m^2-z)(1-u^2)+(\alpha-\underline\alpha)u^2 .
\]
Both differences are nonnegative. Replacing the base by its upper
bound thus increases each integral. Since $Q_0,Q_1\le1$, lowering
the exponent to its integer bound at $n_l$ increases it once more.

\paragraph{Rational bounds on each box.}
Take $S=2^{50}$ and choose all spatial endpoints to be integer
multiples of $1/S$. We round $a_l^2,m_l^2$ down and $a_h^2,m_h^2$
up. Rounding the coefficients of a quadratic upwards preserves an
upper bound for its moment because all three basis functions are
nonnegative on $[0,1]$. In $Q_0$ and $Q_1$ the rounded value
$w=1-\lfloor S m_l^2\rfloor/S$ still lies in $[0,1]$ and increasing $w$ preserves both upper bounds.
The lower bound for $\alpha$ gives
\[
 q^2\le 1-\underline\alpha-a_l^2(1-m_h^2),\qquad
 m^2-\alpha\le m_h^2-\underline\alpha .
\]
We use nonnegative upper bounds for both radicands. Since $0\le q\le1$,
we may replace $n$ by $n_l$ when bounding $q^n$. For the factors
outside the two integrals we use $n_h,m_h,1-m_l^2$ and $C_{n_h-1}$.
The last substitution is valid because $(1+2/u)^{u/2}$ is increasing
for $u>0$, as follows from $\log(1+t)\ge t/(1+t)$. These substitutions
give an upper bound for $q^n+E_0+E_1$. If it lies strictly below
$1-a_hm_h$ rounded down, then~\eqref{finiteineq} fails everywhere
in the box.

For clarity, before outward rounding let $N_h=n_h-1$,
$k_0=\lfloor(n_l-3)/2\rfloor$, $k_1=\lfloor(n_l-1)/2\rfloor$ and put
\[
 u_q=\max\{0,1-\underline\alpha-a_l^2(1-m_h^2)\},\qquad
 u_d=\max\{0,m_h^2-\underline\alpha\}.
\]
Every putative counterexample in the box would then satisfy
\begin{align*}
 1-a_hm_h &\le 1-am\le q^n+E_0+E_1 \\
 &\le u_q^{n_l/2}
   +\frac{n_hN_h}{2} C_{N_h}m_h^4(1-m_l^2)M_{2,k_0}(Q_0) \\
 &\hspace{14mm}{}+2n_hm_h\sqrt{u_d}\,M_{0,k_1}(Q_1).
\end{align*}
All quantities on the last two lines are fixed by the endpoints and
the reflection comparison. Rounding these lines upwards and the first
quantity downwards preserves the implication for every point of the box.

To evaluate~\eqref{eq:moment-sum}, let $A',B',C'$ be the integer
numerators of the rounded coefficients with common denominator $S$.
For $P(z)=(A'+B'z+C'z^2)^k=\sum d_h z^h$, comparing coefficients in
$(A'+B'z+C'z^2)P'=k(B'+2C'z)P$ gives, with $d_{-1}=0$ and
$d_0=(A')^k$,
\[
 h A'd_h=(k-h+1)B'd_{h-1}+(2k-h+2)C'd_{h-2}\quad(1\le h\le2k).
\]
When $A'>0$ every division in this recurrence is exact, even when
the first summand on the right is negative. If $A'=0<C'$, we reverse
both the coefficients and the corresponding beta-integral weights.
If $A'=C'=0$, only the central term $(B')^k$ contributes, and when
$k=0$ we have $P=1$ in all cases. We clear the denominators of the
beta integrals with their least common multiple, account for the factor
$S^{-k}$, and round the resulting moment upwards to a multiple of
$1/S$. Upper bounds involving a square root are obtained in the same
units by computing $\lceil\sqrt K\rceil$ with integer arithmetic for
$K\ge0$. Thus every comparison in the subdivision is an integer
comparison.

The indices in the coefficient recurrence come from taking the coefficient of \(z^{h-1}\) on each side of the differentiated identity. With coefficients outside the range taken as zero, the terms before rearrangement are
\[
 hA'd_h+(h-1)B'd_{h-1}+(h-2)C'd_{h-2}
       =kB'd_{h-1}+2kC'd_{h-2}.
\]
Moving the last two terms on the left to the right gives the recurrence above. For \(A'>0\) the numerator is therefore \(hA'\) times the integer coefficient \(d_h\). If the coefficients are reversed to handle \(A'=0<C'\), the exponent \(h\) is exchanged with \(2k-h\), which explains the corresponding reversal of the integration weights.

We begin with all integer degrees in $[5,256]$ and the entire unit
square for $(a,m)$. Three inequalities suffice to exclude a box. The
first compares the lower bound for $\alpha$ obtained by reflection
with
\[
 m_h^2-\operatorname{dist}(m_h,[a_l,a_h])^2
   =\max_{a_l\le a\le a_h,\ m_l\le m\le m_h}(2am-a^2).
\]
Since $a\ge0$, the maximum over $m$ occurs at $m_h$, and completing
the square gives the stated formula. A lower bound greater than this
maximum excludes the box. The second test uses Corollary~\ref{q-small}
and the monotonicity of $R_N$ to exclude the box whenever the upper
bound for $q^2$ is at most $\min(1/64,R_{n_l-1})^2$. The third
excludes the box when~\eqref{finiteineq} fails throughout it.

If all three tests are inconclusive, put $w=\max(a_h-a_l,m_h-m_l)$.
When $n_h>n_l$ and either $n_h-n_l>n_lw$ or $w\le1/S$, we divide
the degree interval at its midpoint into two consecutive integer
intervals. In all other cases we bisect the longer side of the spatial
rectangle, choosing the $a$-interval when the sides have equal length.
Each subdivision covers the whole preceding box, including the common
boundary when a spatial side is bisected.

The calculation terminates with \MomentLeaves{} boxes at a maximum
depth of \MomentDepth. Of these, \MomentPolar{} are excluded by
reflection and the geometric bound on $\alpha$, \MomentNear{} by
Corollary~\ref{q-small}, and \MomentOrigin{} by~\eqref{finiteineq}.
The program computes the subdivision from the inequalities themselves.
It would report an unresolved box if it reached a single degree and
the minimum spatial width without an exclusion. Every terminal box
satisfies one of the three tests, so induction over the subdivisions
excludes the entire initial domain. This proves the finite-degree
theorem. Degree five also follows from the fourth-moment estimate in
Section~1.

The author's Python script carries out the
calculation. A second implementation uses the same subdivisions and exclusion tests but evaluates the moments in the ordinary power basis by integration by parts. It gives the same three counts and the same maximum depth. The author is willing to supply the scripts when requested.

\paragraph{Completion of the proof.}
The estimates in this section exclude~\eqref{hyp} in degrees five through 256. Together with the fourth-power calculation in Section~1 and Theorem~\ref{large}, this excludes~\eqref{hyp} at every simple interior zero. The endpoint cases give the equality statement, and undoing the rotation and the monic normalization proves Theorem~\ref{main}.

\newcounter{precontrolequations}
\newcounter{precontrolfigures}
\setcounter{precontrolequations}{\value{equation}}
\setcounter{precontrolfigures}{\value{figure}}
\begingroup
\setcounter{equation}{0}\renewcommand{\theequation}{F\arabic{equation}}\renewcommand{\theHequation}{feedback.\arabic{equation}}
\setcounter{figure}{0}\renewcommand{\thefigure}{F\arabic{figure}}\renewcommand{\theHfigure}{feedback.\arabic{figure}}

\section{Pole geometry and feedback realization}\label{sec:feedback}
Suppose that \(p(\xi)=\prod_{j=1}^n(\xi-\lambda_j)\) is monic, \(n\ge2\), and write \(q=p'/n\). The zeros of \(p\) are the poles of a minimal all-pole realization. A characteristic polynomial will also be used when a system has internal modes which disappear from a specified transfer.

The conclusion of Sendov's theorem is local. If all poles lie in the closed unit disk, every marked pole has some critical point within distance one. The estimate needed for the worst pole of a stable system follows from the convex hull inclusion of Gauss and Lucas \cite{totik}. To see the geometry, let \(c\) be a zero of \(p'\) with \(p(c)\ne0\). Conjugating \(0=\sum_j(c-\lambda_j)^{-1}\) gives
\[
 c=\frac{\sum_j\lambda_j |c-\lambda_j|^{-2}}{\sum_j |c-\lambda_j|^{-2}} .
\]
The weights are positive. A multiple zero of \(p\) also occurs in \(p'\), with multiplicity reduced by one. Consequently, on writing \(\rho(p)=\max_j|\lambda_j|\) and \(\alpha(p)=\max_j\operatorname{Re}\lambda_j\), we obtain
\[
 1-\rho(q)\ge1-\rho(p),\qquad -\alpha(q)\ge-\alpha(p).
\]
The first quantity is the radial pole margin for a Schur polynomial. The second is the distance of the rightmost pole of a Hurwitz polynomial from the imaginary axis. For instance, differentiating \((s+1)(s+3)\) puts the remaining pole at \(-2\). The continuous-time margin grows from one to two even though the original pole at \(-3\) was farther from the origin. Each differentiation preserves the centroid, since the coefficient of degree \(n-2\) in \(q\) is \((n-1)/n\) times the coefficient of degree \(n-1\) in \(p\). The last linear factor in the derivative sequence has its zero at this centroid.

\begin{corollary}[The poles near a marked open-loop pole]\label{control:marked}
Let all zeros of the monic polynomial \(p\) lie in the disk \(|z-b|\le d\), with \(d>0\), and let \(q=p'/n\). At every original pole \(\lambda\), the poles \(c_j\) prescribed by \(q\) satisfy
\[
 \frac1{n-1}\sum_{j=1}^{n-1}\frac{1}{|\lambda-c_j|^2}\ge\frac1{d^2}.
\]
Equality holds precisely when \(|\lambda-b|=d\) and \(p(z)=(z-b)^n-(\lambda-b)^n\). In particular one of the derivative poles lies within distance \(d\) of the marked pole. If \(K\) is the compression from~\eqref{app:compression}, then at every simple marked pole
\[
 \| (\lambda I-K)^{-1}\|_2\ge d^{-1},\qquad
 \| (\lambda I-K)^{-1}\|_{\mathrm F}\ge \sqrt{n-1}\,d^{-1}.
\]
\end{corollary}
\begin{proof}
Apply Theorem~\ref{main} to \(d^{-n}p(b+dz)\). Its critical displacements from the transformed marked pole are \((\lambda-c_j)/d\), which gives both the factor \(d^{-2}\) and the equality condition. At a multiple pole one displacement vanishes and the reciprocal sum is infinite. At a simple pole the inverse exists, with eigenvalues \((\lambda-c_j)^{-1}\). Its largest singular value dominates their largest modulus, and in a unitary Schur basis its squared Frobenius norm is at least the sum of their squared moduli. These two observations give the stated matrix estimates.
\end{proof}
This is the local information supplied by the quadratic Tang--Zhang inequality~\cite{z}. It concerns the \(n-1\) poles prescribed by the derivative, independently of the additional pole used to realize them by feedback. The pole margins discussed above concern the extremal modulus and real part.

\subsection{Putting the derivative in a feedback loop}
Use \(D x=\dot x\) in continuous time and \(D x_k=x_{k+1}\) in discrete time. Consider the controllable companion realization
\[
 D x_1=x_2,\quad \ldots,\quad D x_{n-1}=x_n,\qquad
 D x_n=-\sum_{j=0}^{n-1} a_j x_{j+1}+u,
 \qquad p(\xi)=\xi^n+\sum_{j=0}^{n-1}a_j\xi^j .
\]
Fix a real stable pole \(\gamma\) and form the polynomial of degree at most \(n-1\)
\[
 k(\xi)=(\xi-\gamma)q(\xi)-p(\xi)=\sum_{j=0}^{n-1} k_j\xi^j .
\]
This is the companion-coordinate form of classical pole assignment~\cite{wonham}. Here the target is fixed by the derivative and a single additional pole. Apply the static law \(u=v-\sum_j k_j x_{j+1}\). In the last equation each coefficient \(a_j\) becomes \(a_j+k_j\), so the characteristic polynomial of the entire closed loop is exactly \((\xi-\gamma)q(\xi)\). There are still \(n\) internal poles. The extra pole can be chosen at zero in discrete time, or to the left of the rightmost zero of \(q\) in continuous time. In the discrete zero-pole case,
\[
 k_j=-\left(1-\frac jn\right)a_j,\qquad 0\le j<n.
\]
The calculation requires the nominal coefficients and measured companion states, or their equivalents after a known change of coordinates. It calls for no dynamic controller state.

For the measured channel \(y=x_2-\gamma x_1\), the first state equation gives \(Y=(\xi-\gamma)X_1\) with zero initial state. Keeping this output both before and after feedback gives the exact transfer functions
\[
 \frac{Y}{U}=\frac{\xi-\gamma}{p(\xi)},\qquad
 \frac{Y}{V}=\frac{\xi-\gamma}{(\xi-\gamma)q(\xi)}
             =\frac1{q(\xi)} .
\]
Thus the claimed replacement has a physical output interpretation when that linear combination is measured. If the plant output is instead \(x_1\), its closed-loop transfer is \(1/[(\xi-\gamma)q(\xi)]\). Its relative degree is \(n\), while \(1/q\) has relative degree \(n-1\). In this all-pole channel a causal proper reference-to-input path cannot remove that difference. An observer can supply missing states under observability, at the cost of its own error dynamics. The internal model used for robust tracking or disturbance rejection serves a separate purpose and depends on the reference or disturbance to be rejected \cite{francis}.

\paragraph{A two-state implementation.}
Take \(p(z)=(z-\tfrac12)(z+\tfrac34)=z^2+\tfrac14 z-\tfrac38\), so \(q(z)=z+\tfrac18\). With \(\gamma=0\), \(k_0=3/8\) and \(k_1=-1/8\), the gain replaces the original last row \((3/8,-1/4)\) by \((0,-1/8)\). The full pole radius falls from \(3/4\) to \(1/8\). At the same measured output \(y=x_2\), the old transfer was \(z/p(z)\) and the new one is \(1/(z+1/8)\).
\subsection*{Feedback on the mean coordinate}
There is also a spatial interpretation of the same design. Take \(A=\operatorname{diag}(\lambda_1,\ldots,\lambda_n)\), \(e=n^{-1/2}(1,\ldots,1)^\top\), with actuator direction \(e\). Write \(P=I-ee^*\) and choose orthonormal columns \(Q\) spanning \(e^\perp\). The identity \(\det(\xi I-Q^*AQ)=q(\xi)\) is the differentiator construction of Davis and Pereira~\cite{davis,pereira}. Appendix~\ref{app:control} derives it including multiplicities. Now set \(u=(\gamma e^*-e^*A)x+v\). At zero external input the closed-loop matrix is \(M=PA+\gamma ee^*\), with blocks
\[
 [e,Q]^*M[e,Q]=
 \begin{pmatrix}\gamma&0\\Q^*Ae&Q^*AQ\end{pmatrix}.
\]
A trajectory starting with \(e^*x=0\) stays there and has exactly the characteristic polynomial \(q\). The transfer from \(v\) to the mean coordinate \(e^*x\) is \(1/(\xi-\gamma)\). This is why the output needs to be specified along with the characteristic polynomial.

For additive matrix perturbations in the sense of Hinrichsen and Pritchard \cite{hp}, we obtain an additional guarantee. For a normal Schur plant in these coordinates put \(R=\max_j|\lambda_j|\). With \(\gamma=0\), \(\|PA\|_2\le R\). It follows directly that every complex additive matrix error of norm less than \(1-R\) preserves Schur stability, the exact open-loop radius for normal \(A\). There is a strict improvement if every outermost root is simple. Indeed the barycenter formula places every new critical point strictly inside the circle of radius \(R>0\) unless it is an inherited multiple boundary root. If the matrix stability radius of \(M\) equalled \(1-R\), compactness of the unit circle and unit sphere would give \(|z|=1,\|v\|=1\) with
\[
 1-R=\|zv-Mv\|_2\ge1-\|Mv\|_2\ge1-R .
\]
Equality forces \(Mv=Rzv\), contradicting the location of the full closed-loop eigenvalues. A repeated outermost root, on the other hand, persists and realizes equality. The Euclidean norm and normal realization enter the proof explicitly. For a real polynomial each conjugate pair can equivalently be represented by a real planar block, using orthogonal modal coordinates in the real system.
\subsection*{Approaching the derivative poles with one measured output}

There is also a continuous-time experiment in which the controller needs
only one measurement. Let \(A=\operatorname{diag}(r_1,\ldots,r_n)\),
\(p(s)=\det(sI-A)\), and \(e=\mathbf 1/\sqrt n\). Suppose the input acts
along \(e\) and the sensor measures \(y=e^*x\). The feedback \(u=-ky\)
then gives \(A_k=A-kee^*\). Its implementation needs only the measured
output and a positive scalar gain. The matching of the actuator and sensor
directions matters because it gives equal weights to the modes. Indeed,
the determinant lemma and logarithmic differentiation give
\[
 \det(sI-A_k)=p(s)\bigl(1+k e^*(sI-A)^{-1}e\bigr)
            =p(s)+\frac{k}{n}p'(s).
\]
This is the root locus of the transfer \(p^{\prime}/(np)\) in the sense of Evans~\cite{evans}. After division by \(k\), these polynomials converge uniformly on compact
sets to \(p'/n\). On small disjoint circles around its distinct zeros,
Rouch\'e's theorem therefore gives the same number of enclosed poles,
counting multiplicity, for every sufficiently large \(k\). This accounts
for \(n-1\) poles. Since \(\operatorname{tr} A_k=\sum_jr_j-k\) and the
zeros of \(p'\) sum to \((n-1)\sum_jr_j/n\), the remaining pole is
\(-k+\sum_jr_j/n+o(1)\).

If \(\max_j\operatorname{Re}r_j=-d<0\), then
\[
 \operatorname{Re}(x^*A_kx)
 =\sum_j |x_j|^2\operatorname{Re}r_j-k|e^*x|^2
 \le -d\|x\|_2^2 .
\]
Thus every positive gain preserves the decay estimate
\(\|x(t)\|_2\le e^{-dt}\|x(0)\|_2\). In discrete time the same
characteristic polynomial has the same escaping pole, which eventually
leaves the unit disk.
\subsection{Robustness under coefficient perturbations}

Stability radii for polynomial matrices with prescribed coefficient perturbations are treated by Genin, \c{S}tefan and Van Dooren~\cite{genin}. In the scalar companion form one can calculate the radius directly. To obtain a robustness statement, we must also specify how a perturbation enters the system. Let \(P\) be a monic Schur polynomial of degree \(n\), and suppose that its nonleading coefficients undergo additive complex perturbations. Write
\[
 P_\delta(z)=P(z)+\sum_{j=0}^{n-1}\delta_j z^j,\qquad
 d(P)=\min_{|z|=1}|P(z)|.
\]
Denote by \(R_r(P)\) the infimum of \(\|\delta\|_r\) over perturbations for which \(P_\delta\) ceases to be Schur. If \(1\le r\le\infty\) and \(1/r+1/t=1\), then
\[
 R_r(P)=\frac{d(P)}{n^{1/t}},
 \qquad n^{1/\infty}=1 .
\]
Indeed, on the unit circle H\"older's inequality gives
\(\left|\sum_j\delta_jz^j\right|\le n^{1/t}\|\delta\|_r\).
When this bound is less than \(d(P)\), Rouch\'e's theorem places all \(n\) zeros of \(P_\delta\) in the open unit disk. Conversely, choose \(\zeta\) on the circle with \(|P(\zeta)|=d(P)\), and take
\(\delta_j=-P(\zeta)\overline{\zeta}^{\,j}/n\).
At \(\zeta\) the \(n\) contributions have the same phase and sum to \(-P(\zeta)\). This perturbation has norm \(d(P)/n^{1/t}\), so the stated radius is attained.

Now take \(P=p\), set \(q=p'/n\), and return to degree \(n\) by choosing the auxiliary pole at zero. The resulting characteristic polynomial is \(D(z)=zq(z)\). The minimum-modulus inequality of Aziz and Dawood~\cite{aziz} gives \(d(q)\ge d(p)\). Here it also follows directly from the geometry of differentiation. For \(|\lambda|<d(p)\), Rouch\'e's theorem shows that \(p(z)-\lambda z^n\) has \(n\) zeros in the open unit disk. By Gauss--Lucas its derivative \(n(q(z)-\lambda z^{n-1})\) has all its zeros there as well. Thus \(q(\zeta)/\zeta^{n-1}\) lies outside the open disk of radius \(d(p)\) at every \(|\zeta|=1\). Since \(|D(\zeta)|=|q(\zeta)|\), we obtain
\[
 R_r(D)=\frac{d(q)}{n^{1/t}}\ \ge\ \frac{d(p)}{n^{1/t}}=R_r(p).
\]
This comparison applies with the feedback gain fixed at its nominal value. A plant coefficient $a_j+\delta_j$ yields the closed-loop coefficient $a_j+k_j+\delta_j$, so precisely the same perturbation vector acts on both characteristic polynomials. Both radii concern the same number of coefficients and the same norm. In particular, \(r=\infty\) measures the largest common tolerance on the complex coefficients. The calculation gives a sufficient tolerance for real perturbations too, although the perturbation attaining the complex radius can have nonreal coefficients.

The same comparison can be iterated without changing the uncertainty class. If \(q_k=(n-k)!\,p^{(k)}/n!\), then \(q_{k+1}=q_k'/(n-k)\). The degree-\(n\) targets \(D_k(z)=z^kq_k(z)\) satisfy \(d(D_k)=d(q_k)\), so \(R_r(D_k)\) increases weakly at every step. The perturbations in this statement apply to all coefficients below the leading one, including the coefficients that vanish in the nominal target. For an uncertain physical parameter entering several coefficients simultaneously, its effect on those coefficients must be determined before applying this radius.

There is a corresponding input--output interpretation. Take the same disturbance numerator \(B(z)=z\) for each denominator. In the scalar difference equation this choice gives the disturbance a delay of \(n-1\) samples relative to the output. Since \(|B|=1\) on the unit circle, the largest frequency gains are \(\|z/p\|_\infty=1/d(p)\) and \(\|z/D\|_\infty=1/d(q)\). The differentiated design therefore weakly reduces the induced \(\ell^2\) gain in this channel.

A change in the disturbance channel can reverse that comparison. For \(0<a<1\), take \(p(z)=z(z-a)\), \(D(z)=z(z-a/2)\), and use the same numerator \(B(z)=z-1\) in both systems. This numerator takes the first difference of the input, together with a delay. For \(0\le b<1\),
\[
 \left|\frac{e^{i\theta}-1}{e^{i\theta}-b}\right|^2
 =\frac{2(1-\cos\theta)}{(1-b)^2+2b(1-\cos\theta)} .
\]
The fraction increases with \(1-\cos\theta\), so its maximum is \(4/(1+b)^2\). The outer pole radius decreases from \(a\) to \(a/2\), yet the peak gain increases from \(2/(1+a)\) to \(2/(1+a/2)\). The difference already occurs for two stable systems of the same order driven through the same first-difference channel.

In fact, a real constant feedback gain attains the loss of stability in this last example. For the characteristic polynomial \(P+kB\), stability holds whenever \(|k|<1/\|B/P\|_\infty\) by Rouch\'e's theorem. At the positive endpoint \(k=-P(-1)/B(-1)\) a pole reaches \(-1\). Thus the first-difference channel has a smaller feedback gain margin after differentiation, even though the independent coefficient tolerances proved above improve weakly.
\section{Sampling and differentiation in either order}\label{sec:sampling-control}
Let \(\lambda_1,\ldots,\lambda_n\) be the poles of a continuous-time system and let \(T>0\) be the sampling period. Exact sampling of the autonomous dynamics sends each pole \(\lambda_j\) to \(e^{T\lambda_j}\). There are consequently two derivative designs to compare. One differentiates the continuous-time characteristic polynomial and samples its critical points. The other differentiates the polynomial with roots \(e^{T\lambda_j}\). We compare the resulting pole sets by their asymptotic decay rate per unit time, \(-\log\rho/T\), where \(\rho\) is their largest modulus. A set consisting entirely of zero poles has infinite rate in this convention.

For two poles, put \(c=(\lambda_1+\lambda_2)/2\) and \(d=(\lambda_1-\lambda_2)/2\). Differentiation before sampling gives the pole \(e^{Tc}\). Differentiation after sampling gives
\[
 \frac{e^{T\lambda_1}+e^{T\lambda_2}}{2}
     =e^{Tc}\cosh(Td).
\]
Since \( |\cosh(x+iy)|^2=1+\sinh^2x-\sin^2y\), sampling first gives the larger decay rate precisely when
\[
 |\sin(T\operatorname{Im}d)|>|\sinh(T\operatorname{Re} d)|.
\]
For distinct real poles, differentiation before sampling gives the larger rate. For a conjugate pair \(-a\pm i\omega\), with \(a>0\), the two poles to compare are \(e^{-aT}\) and \(e^{-aT}\cos(\omega T)\). Sampling first improves the rate whenever \(\omega T\notin\pi\mathbb Z\). At a quarter-period the latter pole is zero.

For real poles, the non-strict ordering follows from Bourin's compression inequality~\cite[Theorem~1.1]{bourin}. Apply it to the Hermitian diagonal matrix of the poles on the subspace \(e^\perp\), with \(f(x)=e^{Tx}\). The ordered eigenvalues of \(f(Q^*AQ)\) are bounded above by those of \(Q^*f(A)Q\). The differentiator identity identifies these as the two derivative pole sets. For distinct real poles the strict ordering can also be seen in each gap by the following scalar argument. Suppose
\(\lambda_1<\cdots<\lambda_n<0\), and let \(c_j\in(\lambda_j,\lambda_{j+1})\) be the unique critical point in this interval. Let \(\zeta_j(T)\) be the critical point of the sampled polynomial between \(e^{T\lambda_j}\) and \(e^{T\lambda_{j+1}}\). To compare the corresponding poles in the same units, put \(x_j(T)=T^{-1}\log\zeta_j(T)\). We claim that
\[
       c_j<x_j(T)<\lambda_{j+1}.
\]
Indeed, logarithmic differentiation of the two polynomials gives
\[
 F_0(c_j)=0,\qquad F_T(x_j(T))=0,\qquad
 F_0(x)=\sum_k\frac{1}{x-\lambda_k},\quad
 F_T(x)=\sum_k\frac{T}{1-e^{-T(x-\lambda_k)}} .
\]
The second expression follows by multiplying the sampled logarithmic derivative at \(e^{Tx}\) by \(Te^{Tx}\). Both functions decrease strictly from \(+\infty\) to \(-\infty\) on the gap. For every nonzero real \(y\), the inequalities \(y>1-e^{-y}\) and \(y(1-e^{-y})>0\) give \(1/(1-e^{-y})>1/y\). Applied to \(y=T(x-\lambda_k)\), this yields \(F_T(c_j)>F_0(c_j)=0\). Its zero must therefore lie to the right of \(c_j\). In the last gap both derivative designs improve the original decay rate \(-\lambda_n\), and the continuous-time derivative gives the larger improvement.

There is a complementary result for arbitrarily many oscillatory modes with a common decay rate. Suppose \(\lambda_j=-a+i\omega_j\) for every \(j\). Their continuous-time critical points lie on the same vertical line by Gauss--Lucas, so the decay rate remains \(a\). The sampled poles all lie on the circle of radius \(e^{-aT}\). If the sampled phases are pairwise distinct, each critical point lies strictly inside this circle. Indeed, at a critical point \(w\) distinct from the poles \(z_j\), conjugation of the logarithmic derivative gives
\[
 w=\frac{\sum_j z_j|w-z_j|^{-2}}{\sum_j |w-z_j|^{-2}} .
\]
Every weight is positive. A convex combination of distinct points on the circle with these weights has strictly smaller modulus. The sampled derivative thus has decay rate strictly greater than \(a\). If two phases coincide after sampling, the resulting multiple pole persists in the derivative and its modulus keeps the decay rate exactly at \(a\).
\subsection{A damped oscillator in continuous and sampled time}
\label{sec:oscillator}

Consider a damped oscillator with input along its first coordinate,
\[
 A=\begin{pmatrix}-1&-1\\1&-1\end{pmatrix},\qquad b=(1,0)^{\mathsf T}.
\]
We compare each closed loop with its own open loop using the unstructured complex stability radius in the Euclidean operator norm. This is the smallest norm of an additive state-matrix perturbation that places an eigenvalue on the stability boundary. For a stable matrix \(M\) it equals \(\min_{\zeta\in\Gamma}\sigma_{\min}(\zeta I-M)\), where \(\Gamma\) is the imaginary axis in continuous time and the unit circle in discrete time. The formula follows by taking a rank-one perturbation along the singular vectors of \(\zeta I-M\) at its smallest singular value.

In continuous time \(p(s)=(s+1)^2+1\), so \(p'(s)/2=s+1\). Prescribing \(\dot x_1=-2x_1\) requires the static control \(u=-x_1+x_2\). It gives
\[
 A_c=\begin{pmatrix}-2&0\\1&-1\end{pmatrix}.
\]
The rightmost real part remains \(-1\). The normal open-loop matrix has stability radius one, whereas the closed-loop radius drops to \(\sqrt{3-\sqrt5}\simeq0.8740\). Indeed, with \(x=\omega^2\) the smaller squared singular value of \(i\omega I-A_c\) is \(3+x-\sqrt{5+x}\), which increases for \(x\ge0\).

Now sample the same open-loop plant with a zero-order hold and period \(T=\pi/2\). Put \(R=e^{-T}\) and let \(J\) denote counterclockwise rotation by \(\pi/2\). Since \(e^{At}=e^{-t}(I\cos t+J\sin t)\),
\[
 A_d=RJ,\qquad b_d=\int_0^T e^{At}b\,dt
       =\tfrac12(1+R,\,1-R)^{\mathsf T}.
\]
In the orthonormal basis \([e,Je]\), with \(e=b_d/\|b_d\|_2\), the state matrix is still \(RJ\) and the input acts only in the first coordinate. Its characteristic polynomial is \(z^2+R^2\), whose normalized derivative is \(z\). Choose the auxiliary pole \(H=e^{-2T}=R^2\), the sample of the continuous auxiliary pole. For coordinates \((\xi_1,\xi_2)\), the static control \(u=(H\xi_1+R\xi_2)/\|b_d\|_2\) gives
\[
 M_H=\begin{pmatrix} H&0\\ R&0\end{pmatrix}.
\]

The open sampled radius is \(1-R\). For \(M_0\), direct evaluation on the unit circle gives \(r(M_0)=(\sqrt{4+R^2}-R)/2>1-R/2\). Distance to the set of unstable matrices is 1-Lipschitz, so
\[
 r(M_H)\ge r(M_0)-\|M_H-M_0\|_2>1-R/2-R^2>1-R ,
\]
since \(R=e^{-\pi/2}<1/2\). For the actual value, the Gram matrix of \(zI-M_H\) on \(|z|=1\) has diagonal entries \(|z-H|^2+R^2,1\) and off-diagonal entries of modulus \(R\). A unitary change of phase makes these latter entries constant. The smallest eigenvalue is therefore minimized at \(z=1\), where \(|z-H|\) is smallest. This gives \(r(M_H)\simeq0.87774\), compared with \(1-R\simeq0.79212\), an increase of about \(10.8\%\). The continuous radius decreased by about \(12.6\%\). These percentages compare each feedback with its own open-loop perturbation model. The two controllers are designed in their respective time domains.
\subsection*{The same physical actuator error in both loops}
The oscillator lets us go further and vary one dimensionless parameter in precisely the same physical component. Multiply the actuator strength by a real number \(g\), with nominal value one, and hold both nominal controllers fixed. The continuous input matrix becomes \(gb\). Under a zero-order hold the sampled input matrix consequently becomes \(g b_d\). In the coordinates already used, the two closed-loop matrices are
\[
 M_c(g)=\begin{pmatrix}-1-g&g-1\\1&-1\end{pmatrix},\qquad
 M_d(g)=\begin{pmatrix}gH&(g-1)R\\R&0\end{pmatrix},\qquad H=R^2.
\]
The continuous characteristic polynomial is \(s^2+(2+g)s+2\), so the loop is Hurwitz exactly when \(g>-2\). The sampled polynomial is \(z^2-gHz+(1-g)H\). For a monic real quadratic \(z^2+bz+c\), the quadratic Schur--Jury conditions~\cite{jury} are \(1+b+c>0\), \(1-b+c>0\) and \(1-c>0\). These follow, for example, by setting \(z=(1+s)/(1-s)\), multiplying by \((1-s)^2\) and applying the quadratic Hurwitz conditions. Here the three expressions are \(1+H-2gH\), \(1+H\), and \(1-H+gH\), respectively. The exact intervals are therefore
\[
 \text{continuous: } -2<g<\infty,\qquad
 \text{sampled: } 1-1/H<g<(1+H)/(2H).
\]
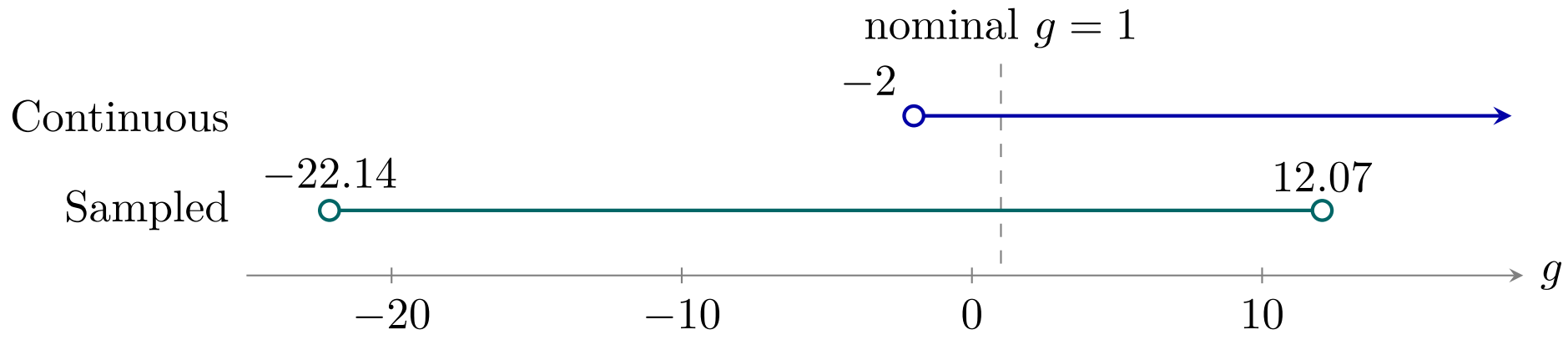
\begin{center}
\begin{tikzpicture}[x=0.23cm,y=0.78cm,>=stealth,font=\small]
 \draw[->,gray] (-25,0)--(19,0) node[right,black] {$g$};
 \foreach \v in {-20,-10,0,10} {
   \draw[gray] (\v,0.08)--(\v,-0.08) node[below,black] {$\v$};
 }
 \draw[gray,dashed] (1,0.12)--(1,2.15);
 \node[above] at (1,2.14) {nominal $g=1$};
 \node[anchor=east] at (-25,1.62) {Continuous};
 \node[anchor=east] at (-25,0.66) {Sampled};
 \draw[->,thick,blue!65!black] (-2,1.62)--(18.6,1.62);
 \draw[thick,blue!65!black,fill=white] (-2,1.62) circle[radius=2.2pt];
 \node[above left] at (-2,1.62) {$-2$};
 \draw[thick,teal!80!black] (-22.14,0.66)--(12.07,0.66);
 \foreach \v in {-22.14,12.07} {
   \draw[thick,teal!80!black,fill=white] (\v,0.66) circle[radius=2.2pt];
   \node[above] at (\v,0.68) {$\v$};
 }
\end{tikzpicture}
\end{center}
With \(H=e^{-\pi}\), the sampled interval is approximately \((-22.14,12.07)\). For a symmetric signed uncertainty \(|g-1|<\varepsilon\), the limiting radii are thus \(3\) in continuous time and \((e^\pi-1)/2\simeq11.07\) after sampling. A different physical question gives the opposite ordering. If the actuator retains its sign and its gain can only increase above one, the continuous controller tolerates every such increase, whereas the sampled loop reaches a pole at \(+1\) when \(g=(1+e^\pi)/2\). The scalar gain and the actuator are identical in these two tests. What changes is the direction of uncertainty for which the designer asks for protection.

\subsection*{Why the two derivative designs require different gains}
One might try to use the continuous feedback unchanged in the sampled implementation. For an oscillator this imposes a restriction even before uncertainty is introduced. Take
\[
 A=-aI+\omega J,\qquad J=\begin{pmatrix}0&-1\\1&0\end{pmatrix},
 \qquad b=(1,0)^{\mathsf T},\qquad a,\omega>0.
\]
Use the same static law \(u=k_1x_1+k_2x_2\) in continuous time and at the sampling instants, holding its value between samples in the latter case. Write \(T>0\) for the sampling period. The normalized derivatives of the open-loop characteristic polynomials have respective zeros \(-a\) and \(c=e^{-aT}\cos(\omega T)\).

\begin{proposition}\label{control:common-gain}
If \(\sin(\omega T)\ne0\), the unique gain assigning both derivative poles is \((k_1,k_2)=(a,\omega)\). Its remaining continuous pole is zero and its remaining sampled pole is one. Thus simultaneous derivative assignment with an unchanged static gain places both systems on their stability boundaries.
\end{proposition}
\begin{proof}
Put \(s=e^{-aT}\sin(\omega T)\) and \(A_T=e^{AT}=cI+sJ\). Write \(B_T=(u_T,v_T)^{\mathsf T}=\int_0^T e^{At}b\,dt\). Then \(AB_T=(A_T-I)b\), whose second coordinate gives \(s=\omega u_T-av_T\). For the row \(k=(k_1,k_2)\), evaluating the two characteristic polynomials at the target poles gives
\[
 \det(-aI-A-bk)=\omega(\omega-k_2),\qquad
 \det(cI-A_T-B_Tk)=s(s+k_1v_T-k_2u_T).
\]
The first expression forces \(k_2=\omega\). The second then becomes \(s(k_1-a)v_T\). To see why \(v_T>0\), write
\[
 v_T=\int_0^T e^{-at}\sin(\omega t)\,dt .
\]
As a function of \(T\) it first increases, and thereafter its local minima occur at \(T=2j\pi/\omega\), \(j\ge1\). At such a minimum its value is \(\omega(1-e^{-aT})/(a^2+\omega^2)>0\). Consequently \(k_1=a\). The first row of \(A+bk\) is now zero, giving the continuous poles \(0,-a\). Finally,
\[
 A_T+B_Tk-I=\left(\int_0^T e^{At}\,dt\right)(A+bk).
\]
Any vector fixed by the continuous flow at this gain therefore also gives a unit eigenvalue of the sampled transition, in agreement with the zero-frequency boundary relation of Ackermann and Hu~\cite{ackhu}. The other sampled eigenvalue is \(c\) by construction.
\end{proof}
When \(\sin(\omega T)=0\), the sampled transition is \(cI\), so the rank-one perturbation \(B_Tk\) then leaves \(c\) as an eigenvalue for every gain.

The same obstruction occurs for two distinct real Hurwitz poles. Write \(A=\alpha I+L\) with eigenvalues \(\alpha\pm d\), \(d>0\). Then \(L^2=d^2I\), \(e^{AT}=cI+sL\) with \(s=e^{\alpha T}\sinh(dT)/d>0\), and \(\Phi=\int_0^T e^{At}\,dt=uI+vL\) with \(v=\int_0^T e^{\alpha t}\sinh(dt)/d\,dt>0\). Comparing the coefficients of \(L\) in \(A\Phi=e^{AT}-I\) gives \(u+\alpha v=s\). For any input column \(b\), the rank-one determinant formula shows that assigning the derivative poles \(\alpha,c\) requires
\[
 kL^{-1}b=-1,\qquad kL^{-1}\Phi b=-s.
\]
Substitution gives \(kb=(u-s)/v=-\alpha\). The continuous trace is then \(\alpha\), so the other continuous pole is zero. The same equilibrium identity as above gives a sampled pole at one.

For the nominal example \(a=\omega=1\), the effect of redesigning the feedback can also be measured directly. Suppose we sample and hold the original continuous gain \((-1,1)\) at \(T=\pi/2\), with actuator multiplier \(g\) and \(R=e^{-\pi/2}\). Substituting the same row in \(RJ+gB_Tk\) gives the characteristic polynomial
\[
 z^2+gRz+R^2-gR.
\]
The three quadratic Schur inequalities give \(1+R^2>0\), \(1+R^2-2gR>0\), and \(1-R^2+gR>0\). Thus the exact interval for the held continuous law is
\[
 R-R^{-1}<g<\tfrac12(R+R^{-1}),\qquad
 (R-R^{-1},(R+R^{-1})/2)\simeq(-4.603,2.509).
\]
For comparison, the interval for the redesigned sampled feedback is obtained by dividing both endpoints by \(R\). Its lower and upper limits are approximately \(-22.14\) and \(12.07\). Measured about \(g=1\), the symmetric tolerance is approximately \(1.51\) for the held continuous law, \(3\) for continuous feedback and \(11.07\) for the sampled derivative design. This third comparison identifies how much of the change can be attributed to redesigning the controller at the chosen sampling period.

\section{A common physical uncertainty model for the two feedback laws}\label{sec:joint}
The singular-value calculation in Section~\ref{sec:sampling-control} compares each matrix with perturbations in its own time domain. To compare the controllers against the same physical change, let the oscillator have damping rate \(a>0\), angular frequency \(\omega>0\), and real actuator multiplier \(g\). The time and gain units are chosen so that the nominal triple is \((1,1,1)\). Here \(\omega\) is the imaginary part of the open-loop pole; the associated second-order scalar polynomial is \(s^2+2as+a^2+\omega^2\). The parameters are constant in each experiment. Only their values change between experiments. Both feedback rows are held fixed throughout.

Let \(b=(1,0)^\mathsf T\), \(J=\bigl(\begin{smallmatrix}0&-1\\1&0\end{smallmatrix}\bigr)\), \(A=-aI+\omega J\) and \(T=\pi/2\). The true transition matrix and the true unit-gain input column over one held interval are
\[
 E=e^{TA},\qquad h=\int_0^T e^{tA}b\,dt=(I-E)(-A)^{-1}b.
\]
The continuous closed-loop generator is \(A+gb f_{\rm c}\), where \(f_{\rm c}=(-1,1)\). The sampled transition is \(M=E+ghf\). Here \(f\) is the nominal sampled row expressed in the original physical coordinates. With \(R=e^{-\pi/2}\), the nominal unit-gain input is \(h_0=(1+R,1-R)^\mathsf T/2\). The coordinate formula of Section~\ref{sec:sampling-control} yields
\begin{equation}
 f=(u,v)=\frac{(R^2h_0+RJh_0)^\mathsf T}{\|h_0\|_2^2},\quad
 u=\frac{R(R^2+2R-1)}{1+R^2},\quad
 v=\frac{R(1+2R-R^2)}{1+R^2}.
 \label{phys:row}
\end{equation}
In particular \(u<0<v\) and \(\|f\|_2=\sqrt2\,R\). Recomputing \(h\) as well as \(E\) matters. A frequency error changes the direction in which a held physical input acts over a sampling interval.

\paragraph{The complete stability region.}
For \(a,\omega>0\) and any real \(g\), the continuous loop is stable exactly when
\[
 2a+g>0,\qquad a^2+\omega^2+g(a-\omega)>0.
\]
For the sampled loop set \(t_0=\operatorname{tr}E\), \(d_0=\det E\), \(t_1=fh\) and \(d_1=f\operatorname{adj}(E)h\). Then the necessary and sufficient conditions are
\begin{equation}
\begin{aligned}
 J_+ &=1-t_0+d_0+g(d_1-t_1)>0,\\
 J_- &=1+t_0+d_0+g(d_1+t_1)>0,\\
 J_0 &=1-d_0-gd_1>0 .
\end{aligned}
\label{phys:jury}
\end{equation}
Indeed the trace of \(M\) is \(t_0+gt_1\) and the rank-one determinant identity gives \(\det M=d_0+gd_1\). Substitution in the three quadratic Schur conditions derived above gives~\eqref{phys:jury}. All three tests are affine in the actuator gain even though their coefficients depend on damping and frequency through an exponential and trigonometric functions. For a fixed pair \((a,\omega)\) they specify the entire interval of stable gains. For a common compact set of such pairs, the exact robust gain interval is their intersection. Stability at sampling instants also gives exponential stability between samples for each fixed plant, because the maps \(e^{\tau A}+g(\int_0^\tau e^{tA}b\,dt)f\) are bounded for \(0\le\tau\le T\).

Ackermann and Hu~\cite{ackhu} developed the comparison in a common physical parameter space. Bernstein and Hollot~\cite{bern-hollot} treated the exponential dependence of the sampled dynamics on plant uncertainty. For a different damped mechanical plant and a digital PD controller, Al-Salem and Fanni~\cite{alsalem} used the Jury test to obtain stability boundaries. Here the two rows come from prescribing the derivative polynomial in their respective time domains. The following calculations locate the first destabilizing plant and gain for these particular rows.

\subsection{An exact radius for simultaneous errors}\label{phys:cube-section}

Let the nominal damping, angular frequency and actuator gain all be one, using the time and gain units of the oscillator above. We hold each feedback row fixed at its nominal value. For \(r\ge0\) consider constant real plant parameters in the common cube
\[
 \mathcal B_r=\{(a,\omega,g):\max(|a-1|,|\omega-1|,|g-1|)\le r\}.
\]
Thus the same relative tolerance in each physical quantity defines the two uncertainty sets. Let \(L=(u+\sqrt2 R)/2\) for the sampled row in~\eqref{phys:row}.
\begin{theorem}\label{thm:joint-cube}
The largest open cubes of stability have respective radii
\[
 \boxed{\quad r_{\rm c}=4\sqrt{2}-5=0.656854249\ldots,\qquad
 r_{\rm d}=\frac{1-L}{1+L}=0.829678401\ldots .\quad}
\]
Every point of the closed cube with a smaller radius is stable. At each stated radius a boundary point has a continuous pole at zero or a sampled pole at one.

\end{theorem}
\begin{proof}
\emph{The common equilibrium condition.}
For any feedback row \((u_0,v_0)\) in the convention \(u_{\rm in}=u_0 x_1+v_0 x_2\), the determinant of the continuous closed-loop generator is
\[
 D=a^2+\omega^2-g(u_0a+v_0\omega).
\]
This same expression controls a sampled pole at one. To see this, represent the plane by \(\mathbb C\) and put
\[
 E=e^{(-a+i\omega)T},\qquad G=\int_0^T e^{(-a+i\omega)t}\,dt
     =\frac{1-E}{a-i\omega}.
\]
The state update is the real-linear map \(x\mapsto Ex+gG\operatorname{Re}((u_0-iv_0)x)\). If \(a>0\), the real matrix of multiplication by \(1-E\) is invertible. The rank-one determinant formula gives
\[
 \det(I-M)=|1-E|^2
   \left(1-g\operatorname{Re}\frac{(u_0-iv_0)G}{1-E}\right)
 =\frac{|1-E|^2}{a^2+\omega^2}D ,
\]
where \(M\) denotes the sampled closed-loop matrix.

For both nominal feedback rows one has \(u_0<0<v_0<2\). Put
\(L_0=(u_0+\sqrt{u_0^2+v_0^2})/2\) and
\(H_0=(\sqrt{u_0^2+v_0^2}-u_0)/2\). Completing the square gives
\[
 D=(\omega-gv_0/2)^2+(a-gL_0)(a+gH_0).
\]
For \(r<1\), the ratio \(a/g\) on the cube is bounded below by \((1-r)/(1+r)\). Thus \(D>0\) throughout if this lower bound exceeds \(L_0\). Where they are equal, choose \(a=1-r\), \(g=1+r\) and \(\omega=gv_0/2\). This point belongs to the cube since \(L_0<v_0/2<1\), the first inequality following by squaring \(\sqrt{u_0^2+v_0^2}<v_0-u_0\). It has \(D=0\). With the continuous row \((-1,1)\), the other Hurwitz condition is \(2a+g>0\) and \(L_0=(\sqrt2-1)/2\), which gives \(r_{\rm c}\). With the sampled row, \(L_0=L\) and the proposed boundary is \(r_{\rm d}\).

\paragraph{Excluding a sampled crossing elsewhere on the circle.}
It remains to check the other two quadratic Schur conditions for all points up to that boundary. We prove them on the slightly larger cube \(\mathcal B_{5/6}\). We use the elementary bounds
\[
 .207<R<.209,\qquad -.11<u<0<v<.3,\qquad
 e^{-\pi/12}<.77,\qquad \cos(\pi/12)<.966 .
\]
For instance \(3.14159<\pi<3.14160\) and the degree-ten exponential series with remainder bounded by its first omitted term divided by \(1-x/12\), for \(0<x<12\), give the exponential bounds. The last bound follows from \(\cos(\pi/12)=(\sqrt6+\sqrt2)/4\). They also give \(L>(\sqrt{2}\cdot.207-.11)/2>1/11\), whence \(r_{\rm d}<5/6\).

Write \(b=|E|\). Throughout the larger cube \(a,\omega,g\in[1/6,11/6]\), so \(0<\omega T<\pi\) and \(b\le e^{-\pi/12}<.77\). Direct expansion of the determinant gives
\[
 \det M=b^2+g\operatorname{Re}((u-iv)\overline{E}G).
\]
Inside this last real part the integrand at time \(t\) is
\(e^{-a(T+t)}[u\cos(\omega(T-t))-v\sin(\omega(T-t))]\).
It is at most \(|u|e^{-a(T+t)}\). Integration therefore yields
\[
 1-\det M\ge (1-b)\left(1+b-\frac{g|u|}{a}b\right)>0 ,
\]
because \(g|u|\le(11/6)\cdot.11<2a\) and \(b<1\).

For the remaining condition let \(Z=G/(1+E)\). The rank-one formula gives
\(\det(I+M)=|1+E|^2[1+g\operatorname{Re}((u-iv)Z)]\). In this range
\[
 \operatorname{Im}Z=
 \frac{\omega\sinh(aT)-a\sin(\omega T)}
 {(a^2+\omega^2)[\cosh(aT)+\cos(\omega T)]}\ge0 ,
\]
using \(\sinh(aT)\ge aT\) and \(\sin(\omega T)\le\omega T\). Since \(v>0\) and \(|G|\le T\), the bracket is at least \(1-g|u|T/|1+E|\). Now \(\cos(\omega T)\ge-\cos(\pi/12)>-.966\), so
\[
 |1+E|^2\ge1+b^2-2(.966)b\ge1+.77^2-2(.966)(.77)>.32^2 .
\]
On the other hand \(g|u|T<(11/6)(.11)(11/7)<.32\). This proves \(\det(I+M)>0\) throughout. All three Schur conditions have now been accounted for, and the equilibrium calculation gives the exact first loss of stability. In particular the sampled loop tolerates a larger simultaneous error even when the three parameters are restricted to positive values. Its first boundary occurs at \((a,\omega,g)\simeq(0.17032,0.25021,1.82968)\).
\end{proof}

The two systems can also be distinguished while the physical actuator stays exactly at its nominal strength. At \((a,\omega,g)=(1/5,1/2,1)\), the continuous determinant is \(-1/100\). The sampled equilibrium factor is positive since \(a^2+\omega^2-u a-v\omega>0.29-0.15>0\). This triple belongs to \(\mathcal B_{5/6}\), on which the preceding proof excludes the other two boundaries. Thus only the sampled loop is stable. In the opposite direction take \(a=\omega=s=1/16\) and \(g=1\). The continuous determinant equals \(2s^2>0\) and the trace is negative. For the sampled loop the factor at \(+1\) has the sign of \(2s-u-v\), which is negative because \(u+v=4R^2/(1+R^2)>1/8\). These two experiments use positive physical parameters and the same nominal actuator strength.

\subsection{Which controller tolerates the larger positive actuator gain?}\label{sec:physical-band}

Suppose that damping and frequency can each vary by the same fraction of their nominal value. For \(0\le\varepsilon\le 1/5\), write \(L=1-\varepsilon\), \(U=1+\varepsilon\) and allow every \((a,\omega)\in[L,U]^2\). The actuator strength \(g\) is now an independent nonnegative parameter. Its upper stability limit is the first gain at which some plant in this square reaches the stability boundary. We vary the physical plant while holding each controller fixed. In particular the sampled input is calculated again for each pair \((a,\omega)\).

With the physical row \(f=(u,v)\) from~\eqref{phys:row}, put \(c=\sqrt2 R\). Define \(\beta=(c-u)/v\) and \(w_*=\min\{U,\beta L\}\).

\begin{theorem}\label{thm:physical-band}
For the continuous and sampled derivative-feedback designs, respectively, the exact upper actuator gains on this family are
\begin{equation}
 G_{\rm c}(\varepsilon)=\frac{1+\varepsilon^2}{\varepsilon},
 \qquad
 G_{\rm d}(\varepsilon)=\frac{L^2+w_*^2}{Lu+w_*v},
 \qquad G_{\rm c}(0)=+\infty .\label{phys-band-ceilings}
\end{equation}
Every gain between zero and the corresponding limit leaves every plant in the square stable. At each finite limit a member reaches the boundary, and that member violates the equilibrium condition at every larger gain. The continuous controller tolerates the larger gain for \(0\le\varepsilon<\varepsilon_c\), while the sampled controller tolerates the larger gain for \(\varepsilon_c<\varepsilon\le 1/5\), where
\[
 \varepsilon_c=\frac{u+v}{2+u-v}=0.1023641\ldots .
\]
The two limits coincide at \(\varepsilon_c\). The worst sampled frequency leaves the upper edge of the square at \(\varepsilon_s=(\beta-1)/(\beta+1)=0.1899758\ldots\). For larger uncertainty it lies strictly inside the frequency interval.
\end{theorem}

\begin{proof}
The continuous closed-loop trace is \(-2a-g<0\). Its determinant is \(a^2+\omega^2+g(a-\omega)\), so a finite upper gain can occur only for \(\omega>a\). Then \((a^2+\omega^2)/(\omega-a)\) increases with \(a\) and decreases with \(\omega\) as long as \(1<\omega/a<1+\sqrt2\). The derivative with respect to \(\omega\) has the sign of \(\omega^2-2a\omega-a^2\). Since \(U/L\le3/2\), the worst plant has \(a=L,\omega=U\), giving the first expression in~\eqref{phys-band-ceilings}. At \(\varepsilon=0\) the determinant stays equal to two.

For the sampled plant set \(r=e^{-a\pi/2}\), \(\theta=\pi\omega/2\), \(S=a^2+\omega^2\), and
\[
 P=au+\omega v,\qquad V=\omega u-av.
\]
In the whole square \(P>0\) and \(V<0\). The Schur condition at \(+1\) is \(gP/S<1\). For fixed \(\omega\), the ratio \(P/S\) decreases with \(a\), since its numerator is positive and has negative derivative in \(a\). At \(a=L\), its derivative with respect to \(\omega\) has the sign of \(v(L^2-\omega^2)-2Lu\omega\), whose unique positive zero is \(\omega=\beta L\). As \(\beta>1\), its constrained maximum occurs at \(w_*\). This gives the candidate sampled limit in~\eqref{phys-band-ceilings}. The nominal plant belongs to every square and already reaches \(+1\) at \((1+e^\pi)/2<13\).

It remains to check that neither of the other two Schur conditions sets an earlier boundary. We do this uniformly for all \((a,\omega)\in[4/5,6/5]^2\) and \(0\le g\le13\). The determinant has constant term \(r^2\) at \(g=0\), and its coefficient of \(g\) is
\[
 \frac rS\bigl[(\cos\theta-r)P+\sin\theta\,V\bigr].
\]
The condition at \(-1\), after division by the positive number \(1+2r\cos\theta+r^2\), is \(1+gH_->0\), where
\[
 H_-=\frac{(1-r^2)P+2r\sin\theta\,V}
              {(1+2r\cos\theta+r^2)S}.
\]
Both expressions follow by writing the held input in complex coordinates as
\((1-re^{i\theta})(a+i\omega)/S\) and resolving it along the row \((u,v)\).

For the signs it suffices to use \(0.10<-u<0.11\), \(0.27<v<0.28\), \(0.29<c<0.30\), all obtained from the displayed nominal row. On the frequency interval, \(|\cos\theta|<1/3\) and \(\sin\theta>9/10\). Since \(P\le (6/5)0.28\) and \(|V|\ge (4/5)0.27\), the determinant coefficient is negative. Consequently the determinant stays below one.

Write \(N_-=(1-r^2)P+2r\sin\theta\,V\). Since \(0<\sin\theta\le1\),
\[
 N_-\ge \omega\bigl[(1-r^2)v-2r|u|\bigr]
               -a|u|-2arv .
\]
For \(4/5\le a\le1\) one has \(r<2/7\) and \(ar<8/35\), since \(ae^{-a\pi/2}\) decreases on this interval. Substitution gives \(N_->0.18\omega-0.238\). For \(1\le a\le6/5\), the bounds \(r,ar<0.21\) give \(N_->0.21\omega-0.25\), which is at least \(0.18\omega-0.238\) in our range. When \(\omega\le1\), the denominator of \(H_-\) is at least \(32/25\) and the numerator exceeds \(-0.094\). When \(\omega\ge1\), the denominator exceeds \((1-4/21)(41/25)>1\) and the numerator exceeds \(-0.058\). Hence \(H_->-1/13\) throughout. The first loss of stability is therefore exactly the one at \(+1\) found above.

While \(w_*=U\), the two gain limits share the numerator \(L^2+U^2\). Their denominators are \(2\varepsilon\) and \(Lu+Uv\), which agree precisely at \(\varepsilon_c\). The latter denominator is larger at zero. Its slope \(v-u\) is less than two, giving the stated reversal. Here \(\varepsilon_c<0.15<\varepsilon_s\), as can also be checked from the elementary bounds above. For \(\varepsilon\ge\varepsilon_s\), \(G_{\rm c}<7\), while \(G_{\rm d}=2L/(c+u)>8\). Thus no second reversal occurs.
\end{proof}

At five percent variation in both damping and frequency the continuous gain limit is \(20.05\), against approximately \(10.852\) for the sampled design. At twenty percent the continuous limit is \(5.2\), whereas the sampled controller tolerates gains up to approximately \(8.594\). In the latter case the sampled instability first occurs at \(a=0.8\) and \(\omega=1.17525\ldots\). Testing only the four corners of the frequency and damping square would miss that plant.

Figure~\ref{fig:physical-gain-crossing} shows where the upper gain limits exchange order.

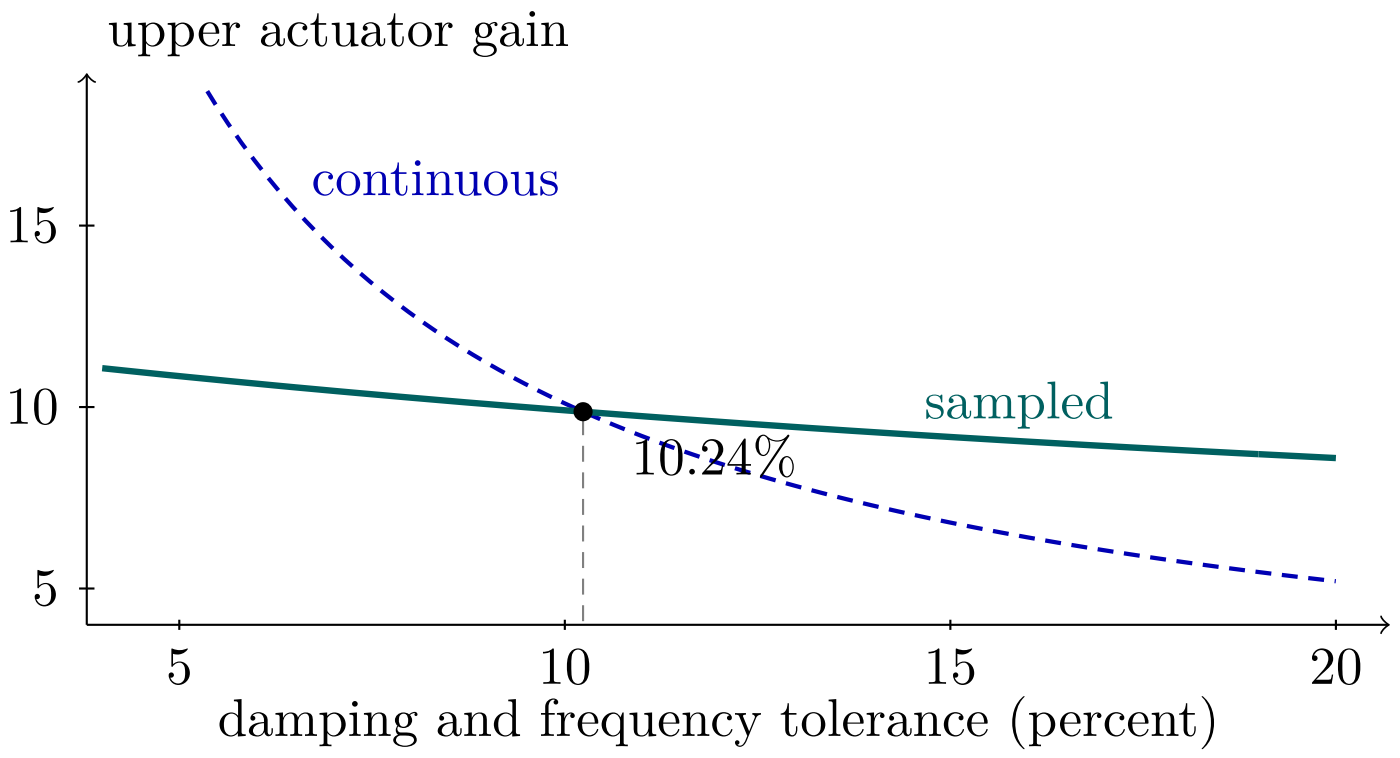
\begin{figure}[htbp]
\centering
\begin{tikzpicture}[x=0.51cm,y=0.24cm,font=\small]
  \begin{scope}
    \clip (3.8,4) rectangle (20.4,18.9);
    \draw[blue!70!black,thick,densely dashed,domain=4.3:20,samples=95]
      plot (\x,{100/\x + \x/100});
    \draw[teal!75!black,very thick,domain=4.0:18.9975825,samples=55]
      plot (\x,{(2+2*(\x/100)^2)/
        ((1-\x/100)*(-0.10780957)+(1+\x/100)*(0.27350490))});
    \draw[teal!75!black,very thick,domain=18.9975825:20,samples=8]
      plot (\x,{2*(1-\x/100)/0.186176542});
  \end{scope}
  \draw[->] (3.8,4)--(20.7,4);
  \draw[->] (3.8,4)--(3.8,19.2);
  \foreach \x in {5,10,15,20}
     \draw (\x,4.14)--(\x,3.86) node[below] {$\x$};
  \foreach \y in {5,10,15}
     \draw (3.9,\y)--(3.7,\y) node[left] {$\y$};
  \draw[gray,densely dashed] (10.2364,4.1)--(10.2364,9.8714);
  \fill (10.2364,9.8714) circle (1.8pt);
  \node[anchor=west,blue!70!black] at (6.45,16.3) {continuous};
  \node[anchor=west,teal!75!black] at (14.4,10.05) {sampled};
  \node[anchor=west] at (10.6,8.65) {\(10.24\%\)};
  \node[anchor=south west] at (3.8,19.1) {upper actuator gain};
  \node[anchor=north] at (12,2.6) {damping and frequency tolerance (percent)};
\end{tikzpicture}
\caption{The continuous (dashed) and sampled (solid) upper gain limits for all plants whose damping and frequency each vary by at most \(\varepsilon\) from one. The two limits meet at \(\varepsilon_c\), and both are computed with the controllers held fixed.}\label{fig:physical-gain-crossing}
\end{figure}

\subsection{An upward drift in the natural rotation frequency}\label{sec:freqdrift}
The two tolerances can also be chosen independently. As an example, suppose \(L=0.99\le a\le1.01\) and \(1\le\omega\le W\le5/2\). Write \(c=\sqrt2 R\) and use the row \((u,v)\) of~\eqref{phys:row}.
\begin{corollary}\label{ctrl:drift}
Put \(w_{\rm c}=\min\{W,L(1+\sqrt2)\}\) and \(w_{\rm d}=\min\{W,L(c-u)/v\}\). The first nonnegative actuator gains for which some member of this family loses stability are
\[
 G_{\rm c}=\frac{L^2+w_{\rm c}^2}{w_{\rm c}-L},\qquad
 G_{\rm d}=\frac{L^2+w_{\rm d}^2}{Lu+w_{\rm d}v}.
\]
They coincide precisely at \(W=L(1+u)/(1-v)=1.215794\ldots\). The continuous gain limit is larger below this frequency bound, and the sampled limit is larger above it.
\end{corollary}
\begin{proof}
The stationary frequencies of the two equilibrium tests are those in the proof of Theorem~\ref{thm:physical-band}, now restricted to \([1,W]\). For the sampled test, \(h_*=(au+\omega v)/(a^2+\omega^2)\) stays positive on the rectangle. Since its numerator decreases with \(a\) and its denominator increases, its maximum lies at \(a=L\). The same endpoint minimizes the finite continuous gain bound. This yields the proposed expressions. To verify that all other conditions in~\eqref{phys:jury} stay strict up to the sampled limit, put \(w=(-A)^{-1}b\), so that \(h=(I-E)w\) and \(h_*=fw\). On the wider frequency interval it is convenient to use
\[
 \|E\|_2<0.22,\quad \|w\|_2<0.72,\quad \|f\|_2<0.30,\quad
 -0.108<u<-0.107,\quad 0.273<v<0.274 .
\]
These follow from~\eqref{phys:row} and \(0.20787<R<0.20789\) by the same exponential-series estimate as in Theorem~\ref{thm:joint-cube}. In particular the sampled upper gain cannot exceed its value at \(W=1\), which is less than 12. Throughout the rectangle \(h_*>1/20\), since
\[
 au+\omega v-\tfrac1{20}(a^2+\omega^2)
 > -0.11(1.01)+0.27\omega-\tfrac1{20}(1.01^2+\omega^2)>0 .
\]
The last quadratic is concave and is positive at both endpoints. From \(f(I+E)^{-1}h=h_*-2f(I+E)^{-1}Ew\), the rank-one determinant formula gives
\[
 \frac{\det(I+M)}{\det(I+E)}
 > 1+g\left(0.05-\frac{2(0.30)(0.22)(0.72)}{0.78}\right)>0
 \qquad (0<g\le12).
\]
Also \(\|\operatorname{adj}(E)\|_2=\|E\|_2\) for this scaled rotation, so
\(\det M<0.22^2+12(0.30)(0.22)(1.22)(0.72)<0.75\).
The case \(g=0\) is stable. Thus the equilibrium test gives the sharp bound over the whole rectangle. Until the first stationary frequency, the gain numerators coincide and their denominators agree exactly at \(W-L=Lu+Wv\). This occurs before the sampled stationary frequency \(1.45437\ldots\). Beyond that frequency the sampled limit stays constant, while the continuous limit continues decreasing to its minimum at \(2.39007\ldots\). 
\end{proof}
For \(W=5/2\) the respective gain limits are \(4.78014\ldots\) and \(10.63507\ldots\). A positive physical gain of five is therefore safe for every sampled plant in this uncertainty set. At \(a=1,\omega=5/2,g=5\) the continuous determinant equals \(-1/4\), and the continuous loop already has a positive real pole.

\subsection*{When does sampled feedback give the larger stability margin?}

The damped oscillator gives an example in which a sampled feedback law has
a larger robust stability margin than its continuous-time counterpart. Both
laws are designed for the same nominal plant and act through the same
actuator. Suppose the damping, natural frequency and actuator gain can each
deviate independently from their nominal values by a constant relative
error. Theorem~\ref{thm:joint-cube} gives the limiting common tolerance as
\(65.6854\ldots\%\) for continuous feedback and \(82.9678\ldots\%\) for
sampled feedback, a difference of about \(17.28\) percentage points.
At a tolerance of \(75\%\), for instance, every plant in the uncertainty
cube is stable under the sampled law. The same cube contains the plant
\((a,\omega,g)=(1/4,1,7/4)\). Its frequency is exactly nominal, yet the
continuous closed-loop determinant is
\(a^2+\omega^2-g(\omega-a)=-1/4\), so one pole is real and positive.

The choice between the two laws also depends on which physical parameter
is expected to vary. Suppose the damping stays within \(1\%\) of nominal
and the normalized frequency lies in \([1,W]\), with \(1\le W\le 5/2\).
If we now ask how far the actuator gain can be increased while preserving
stability for every plant in that range, Corollary~\ref{ctrl:drift} places
the crossover at \(W=1.215794\ldots\). The continuous law admits the
larger gain below this threshold, and the sampled law admits the larger
gain above it. Thus even for the same oscillator, the frequency range
expected in operation can determine which feedback law gives the greater
margin against an increase in actuator gain.

\paragraph{Scope of the comparison.}
For a specified static uncertainty set, the inequalities~\eqref{phys:jury} select the stable plants point by point. Theorem~\ref{thm:joint-cube} measures distance to the first unstable physical plant when all three tolerances grow together. Theorem~\ref{thm:physical-band} holds the plant tolerances fixed and asks how far the actuator gain can be increased. The two questions select different worst plants. Their answers refer to state feedback with the indicated rows, a fixed sampling period and a zero-order hold. Time-dependent parameter changes, disturbance attenuation between samples, and robustness of an observer-based implementation would require corresponding dynamical models.

\clearpage
\endgroup
\setcounter{equation}{\value{precontrolequations}}
\setcounter{figure}{\value{precontrolfigures}}

\section{Conclusion}\label{sec:conclusion}

What first drew me to Sendov's question was the possibility of recognizing a
counterexample by the behaviour of a discrete-time system in which the discre-times system has better stability margins than the continuous one, because the criterions for stability are fundamentally different. For the continuous system it is about having the roots as far as possible in the left half-plane for, while for the discrete-time system it is about having the roots as close as possible to the origin. In this particular context we use the idea that if every critical
point were farther than one from a given simple zero, every eigenvalue of the
reciprocal matrix would lie inside the open unit disk. Its trace gives us a
way to see why this cannot happen. The mean of the reciprocal displacements
determines a path of integration, and their variance bounds the error made
when we replace them by that mean. Together with reflection in the unit
circle and the moment bounds for the finite range of degrees, this gives
another proof of the quadratic inequality of Tang and Zhang, from which
Sendov's distance bound follows.

The same derivative that enters the proof as the numerator of
\(p'/(np)\) can also be used to design a feedback law. For the measured
output in Section~\ref{sec:feedback}, the closed-loop transfer function has
denominator \(p'/n\), while the additional pole of the full state equation
remains available to the designer. The oscillator shows what this freedom
can mean for a physical plant. With the two feedback laws fixed at their
nominal settings, the sampled law tolerates a larger simultaneous error in
damping, frequency and actuator gain. When we ask how far the gain can be
increased, the answer depends on how much frequency drift we allow. The
continuous law tolerates the larger gain over a narrow frequency range.
As the upper frequency bound passes \(1.215794\ldots\) times the nominal
frequency, the ordering reverses (Corollary~\ref{ctrl:drift}).

This suggests choosing the sampling period with the expected variation of the
plant in mind. For a system with several oscillatory modes, can we choose the
period and the additional pole so that sampled feedback tolerates a larger
range of physical parameters than a continuous derivative-based law designed
to meet the same nominal decay requirement? In such a comparison it would
also be worth following the response to a disturbance between sampling
instants, while the actuator holds the value chosen at the last sample.

On the polynomial side, the equality case of Theorem~\ref{main} suggests a
quantitative question. Suppose the zeros of a polynomial of degree \(n\) lie
in the closed unit disk and the mean squared reciprocal distance at a simple
marked zero is at most \(1+\varepsilon\). How close must the zeros be to
the vertices of a regular \(n\)-gon on the unit circle? For fixed \(n\),
compactness and the equality case already imply that, after a suitable
rotation and matching, the largest distance between a zero and its vertex
tends to zero with \(\varepsilon\). What one would like to obtain from the
proof is an explicit rate and to understand how it depends on the degree.
The variance measures how tightly the reciprocal displacements cluster. The
further step would be to turn estimates for their mean and variance into an
estimate for the positions of the original zeros.

\appendix

\section{From polynomial zeros to a discrete-time realization}\label{app:control}

The link with discrete-time control is already visible in the logarithmic derivative of a polynomial. If $p(z)=\prod_{i=1}^n(z-r_i)$ is monic, each root contributes one term to
\begin{equation}
 H(z):=\frac{p'(z)}{np(z)}
 =\frac1n\sum_{i=1}^n\frac1{z-r_i}
 =e^*(zI-A)^{-1}e,
 \qquad A=\operatorname{diag}(r_1,\ldots,r_n),\quad
 e=n^{-1/2}\boldsymbol1.\label{app:transfer}
\end{equation}
Thus $H$ is the transfer function of the discrete-time realization
\[
 x_{k+1}=Ax_k+e u_k,\qquad o_k=e^*x_k.
\]
The diagonal entries $r_i$ are the modes of this realization and lie in the closed unit disk by hypothesis. Its transfer function has zeros at the critical points that survive cancellation in $p'/p$. To recover the full multiplicity of a critical point that is also a root of $p$, we return to the determinant before reducing the rational function.

Differentiation has a direct interpretation as a compression of the diagonal matrix $A$. Let the columns of an $n\times(n-1)$ matrix $Q$ form an orthonormal basis of $e^\perp$ and write $K=Q^*AQ$. In the unitary basis $[e,Q]$, the cofactor formula for the upper left entry of $(zI-A)^{-1}$ gives
\begin{equation}
 \det(zI-K)=p(z)e^*(zI-A)^{-1}e=\frac{p'(z)}n.\label{app:compression}
\end{equation}

To see which minor occurs here, set \(U=[e,Q]\). At a point where \(p(z)\ne0\), write the change of basis explicitly as
\[
 U^*(zI-A)U=
 \begin{pmatrix}
    z-e^*Ae & -e^*AQ\\
    -Q^*Ae & zI_{n-1}-K
 \end{pmatrix}.
\]
Its determinant is still \(p(z)\). Deleting the first row and column gives precisely \(zI_{n-1}-K\), and the cofactor has positive sign. The upper left entry of the inverse is therefore \(\det(zI_{n-1}-K)/p(z)\). Reading this same entry in the original coordinates gives \(e^*(zI-A)^{-1}e\). This identifies the minor in~\eqref{app:compression} without any simplicity assumption on the zeros of \(p\).

The calculation holds away from the roots of $p$ and extends to every $z$ by polynomial identity. Thus the eigenvalues of $K$ are the critical points of $p$ with their algebraic multiplicities. If $Kv=cv$ and $\|v\|=1$, then
\[
 c=v^*Kv=(Qv)^*A(Qv)
   =\sum_{i=1}^n |(Qv)_i|^2r_i,
 \qquad \sum_{i=1}^n |(Qv)_i|^2=1.
\]
This convex combination proves the Gauss--Lucas inclusion. Its weights depend on the eigenvector, so the compression locates each critical point relative to the convex hull of all the roots. To control its distance from a particular root, we must keep that root distinguished in the construction.

Accordingly, fix a root $a$ and write
\[
 p(z)=(z-a)\prod_{k=1}^N(z-b_k),\qquad N=n-1,
 \qquad D=\operatorname{diag}(a-b_1,\ldots,a-b_N),
 \quad J=\boldsymbol1\boldsymbol1^*.
\]
When the marked root is simple, $D$ is invertible and the reciprocal realization takes the form
\begin{equation}
 R_a=(I+J)D^{-1},\qquad
 R_a^{-1}=D\left(I-\frac Jn\right),\qquad
 B_a=aI-R_a^{-1}=\operatorname{diag}(b_k)+\frac{DJ}{n}.
 \label{app:inverse}
\end{equation}
Since $J^2=NJ$, we have $(I+J)^{-1}=I-J/n$. The transpose of $B_a$ is the $D$-companion matrix discussed in the introduction. To obtain its characteristic polynomial, put $h(z)=\prod_k(z-b_k)$ and apply the rank-one determinant formula. This gives
\[
 \det(zI-B_a)
 =h(z)\left(1-\frac1n\sum_k\frac{a-b_k}{z-b_k}\right)
 =\frac{h(z)+(z-a)h'(z)}n
 =\frac{p'(z)}n.
\]

In the middle equality, decompose \(a-b_k=(z-b_k)+(a-z)\). The constant terms sum to \(N=n-1\), and the remaining sum is the logarithmic derivative of \(h\), giving
\[
 1-\frac1n\sum_k\frac{a-b_k}{z-b_k}
 =1-\frac{N}{n}+\frac{z-a}{n}\sum_k\frac{1}{z-b_k}
 =\frac{1}{n}+\frac{z-a}{n}\frac{h'(z)}{h(z)} .
\]

The first expression agrees with the polynomial on the right wherever its denominators are nonzero, and the identity extends to all $z$. It follows from~\eqref{app:inverse} that the eigenvalues of $R_a$ are $y_j=(a-c_j)^{-1}$ with their algebraic multiplicities.

We may also obtain the determinant used in the proof directly from $R_a$. With $g(w)=\prod_k(1-w/(a-b_k))$, the rank-one formula gives, away from the zeros of $g$,
\[
 \det(I-wR_a)
 =g(w)\left(1-w\sum_k\frac1{a-b_k-w}\right)
 =g(w)+wg'(w).
\]
Since $p(a-w)=-wp'(a)g(w)$, differentiation with respect to $w$ now gives
\begin{equation}
 \det(I-wR_a)=\frac{p'(a-w)}{p'(a)}
             =\prod_{j=1}^N(1-wy_j).\label{app:determinant}
\end{equation}
The determinant argument applies equally when $R_a$ has nontrivial Jordan blocks.

For the discrete-time system $u_{k+1}=R_a u_k$, asymptotic stability is equivalent to $\rho(R_a)<1$. In the reciprocal coordinates this is precisely the condition that Sendov's bound fails at $a$, since
\[
 |a-c_j|>1\text{ for every }j
 \quad\Longleftrightarrow\quad
 |y_j|<1\text{ for every }j
 \quad\Longleftrightarrow\quad \rho(R_a)<1.
\]
A critical point at distance strictly less than one corresponds to $\rho(R_a)>1$. When $\rho(R_a)=1$, all trajectories are bounded precisely when every Jordan block belonging to an eigenvalue on the unit circle has size one. The quadratic result proved here follows by excluding~\eqref{hyp} at every interior marked root.

The trace criterion for an isolated root illustrates how this stability question first arose. Suppose that after rotation $a>0$ and that $|a-b_k|>1$ for each of the other roots. Since $|b_k|\le1$, we obtain $\Rea b_k<a/2$. The point $a$ is therefore the unique point of largest real part in the convex hull of the roots. By Gauss--Lucas and $p'(a)\ne0$, every critical point satisfies $\Rea c_j<a$, and hence $\Rea y_j>0$. If $R_a$ were Schur stable, then
\[
 \Rea(1/y_j)=\frac{\Rea y_j}{|y_j|^2}>\Rea y_j
 \quad\text{for every }j.
\]
Summing would yield $\Rea\operatorname{tr}R_a^{-1}>\Rea\operatorname{tr}R_a$. Under the isolation assumption, the reverse weak inequality therefore suffices for Sendov's bound at $a$. For general root configurations we retain the reciprocal matrix and use the quadratic mean of its eigenvalue moduli together with the integral of its determinant.

Consider, for example, $p(z)=(z-a)(z-b)$ with distinct roots in the closed unit disk. Its critical point is $c=(a+b)/2$, and $R_a$ reduces to the scalar $2/(a-b)$. Thus
\[
 |a-c|=\frac{|a-b|}{2}\le1,\qquad
 |R_a|^2=\frac4{|a-b|^2}\ge1.
\]
Taking $a=1/2$ and $b=-1/2$ gives a distance of $1/2$ and a reciprocal eigenvalue of $2$. Taking $a=1$ and $b=-1$ gives equality in both bounds. In higher degree the reciprocal relation still holds for each critical point, while the trace and variance describe how their contributions combine in the integral.

A cubic already shows why the passage from the transfer function to the determinant matters. Take \(p(z)=(z-\tfrac12)(z+\tfrac12)^2\) and mark \(a=\tfrac12\). Direct differentiation and cancellation give
\[
 p'(z)=(z+\tfrac12)(3z-\tfrac12),\qquad
 \frac{p'(z)}{3p(z)}
    =\frac{3z-\tfrac12}{3(z-\tfrac12)(z+\tfrac12)} .
\]
The critical point \(c=-\tfrac12\) is a pole of the reduced transfer function, so looking only for its zeros would recover just \(c=\tfrac16\). In the reciprocal realization, both other roots equal \(-\tfrac12\), giving \(D=I_2\) and
\[
 R_a=\begin{pmatrix}2&1\\1&2\end{pmatrix},\qquad
 R_a\binom{1}{-1}=\binom{1}{-1},\qquad
 R_a\binom{1}{1}=3\binom{1}{1}.
\]
The two eigenvalues recover the two critical displacements separately. They are \(1/(a+\tfrac12)=1\) and \(1/(a-\tfrac16)=3\), so their mean square is \((1+9)/2=5\). One critical point is exactly one unit from the marked zero and the other is only one third of a unit away.

\section{Reflection and the two polynomial evaluations}\label{app:reflection}

For $0<a<1$, evaluating $p$ at the origin and at the reflected point $1/a$ gives estimates in opposite directions. The first evaluation supplies an upper bound and the second a lower bound for integrals of the derivative. In both, the critical points enter through the same determinant~\eqref{app:determinant}.

At the origin, the factorization of the monic polynomial gives $|p(0)|\le a$. At the marked root we have
\[
 |p'(a)|=n\prod_j|a-c_j|
          =\frac n{\prod_j|y_j|}\ge n
\]
by the arithmetic--geometric mean inequality and~\eqref{hyp}. Integrating from the marked root to the origin and then rescaling the interval gives
\[
 -p(0)=p'(a)\int_0^a\prod_j(1-ty_j)\,dt,
 \qquad
 \left|n\int_0^1\prod_j(1-at y_j)\,dt\right|\le1.
\]
Thus the locations of the roots bound the modulus of the complex integral before we estimate the error introduced by replacing the reciprocal critical displacements by their mean.

For the evaluation at $1/a$, put $\beta=1-a^2$. Each root $b$ in the closed unit disk satisfies
\begin{equation}
 |1-ab|^2-|a-b|^2=\beta(1-|b|^2)\ge0.
 \label{app:reflection-factor}
\end{equation}
It follows that each quotient $(1-ab_k)/(a-b_k)$ has modulus at least one. We can evaluate $p(1/a)/p'(a)$ either by factoring $p$ or by substituting $z=a+\beta t/a$ in $\int_a^{1/a}p'(z)\,dz$. The two factorizations give
\begin{align*}
 \frac{p(1/a)}{p'(a)}
 &=\frac\beta{a^n}\prod_k\frac{1-ab_k}{a-b_k},\\
 \frac1{p'(a)}\int_a^{1/a}p'(z)\,dz
 &=\frac\beta{a^n}\int_0^1\prod_j(a+\beta t y_j)\,dt.
\end{align*}
Since $p(a)=0$, the left sides agree. Taking absolute values in the resulting identity and applying the arithmetic--geometric mean inequality to the squared moduli of the derivative factors gives
\[
 \left|\prod_j(a+\beta t y_j)\right|
 \le\left(a^2+2a\beta t\Rea s
                +\beta^2t^2\frac1N\sum_j|y_j|^2\right)^{N/2}.
\]
Hypothesis~\eqref{hyp} bounds the final average by one and yields~\eqref{polar}. The lower bound in that inequality comes from the roots $b_k$, while the upper estimate on the integrand uses the mean square of the reciprocal critical displacements.

The parameter $\alpha=2a\Rea s-a^2$ links this reflected inequality to the origin comparison. Its nonnegativity follows from Gauss--Lucas and the identity
\[
 \alpha=(1-a^2)\left(1-\frac1N\sum_j|y_j|^2\right)
       +\frac1N\sum_j(1-|c_j|^2)|y_j|^2.
\]
The identity $\alpha=|s|^2-|a-s|^2$ gives the other bound, so $0\le\alpha\le|s|^2$ under~\eqref{hyp}. The degree enters through $v=N\alpha$, which measures the effect of this geometric bound on the product of $N$ factors.

\section{Cancellation and the variance of the reciprocal displacements}\label{app:variance}

The normalized trace $s=N^{-1}\operatorname{tr}R_a$ is the mean of the reciprocal critical displacements. Writing $y_j=s+d_j$, with $\sum_jd_j=0$ and $V=N^{-1}\sum_j|d_j|^2$, we obtain
\begin{equation}
 \frac1N\sum_j|y_j|^2=|s|^2+V.\label{app:variance-decomposition}
\end{equation}
It follows from~\eqref{hyp} that $|s|^2+V\le1$. Both terms depend on the eigenvalues of $R_a$. The squared Frobenius norm is the sum of the squared singular values and can exceed $\sum_j|y_j|^2$ when $R_a$ is nonnormal.

Replacing each $y_j$ by $s$ gives the product $(1-ts)^N$, whose integral is explicit. To estimate the change in the integral, interpolate between $s$ and $y_j$ by setting $z_j(h)=s+hd_j$ for $0\le h\le1$, and write
\[
 F(h,t)=\prod_j(1-tz_j(h)),\qquad
 P_j=\prod_{\ell\ne j}(1-tz_\ell(h)),\quad
 P_{jk}=\prod_{\ell\ne j,k}(1-tz_\ell(h)).
\]
Differentiating gives $\partial_hF=-t\sum_jd_jP_j$. The relation $\sum_jd_j=0$ allows us to write the sum in a form that retains the cancellation when the factors are estimated. Indeed,
\[
 \sum_jd_jP_j=\frac1{2N}\sum_{j,k}(d_j-d_k)(P_j-P_k),
 \qquad P_j-P_k=ht(d_j-d_k)P_{jk}.
\]
Together these identities give (V2), in which the deviations enter quadratically. We may now take absolute values and use
\[
 \sum_{j,k}|d_j-d_k|^2=2N^2V.
\]
to bound the error in terms of the variance. For a pair of points $s+d,s-d$, the same cancellation reads
\[
 (1-t(s+d))(1-t(s-d))=(1-ts)^2-t^2d^2.
\]
The terms linear in $d$ cancel. The remaining term contains the complex square $d^2$, whose modulus gives the factor $|d|^2$ in the estimate.

The product over the remaining factors is defined even if some $1-tz_j(h)$ vanishes. An estimate for their mean squared modulus still applies when an individual $|y_j|$ exceeds one. On the real interval this gives (V1) near the concentrated configuration $s=1$. The same cancellation applies for complex $t$ along the contour used below.

\section{Why the contour passes through the conjugate trace}\label{app:contour}

To choose the path of integration, consider the estimate on the mean square of the factors $1-at z_j(h)$. Centering gives
\[
 \frac1N\sum_j|1-at z_j(h)|^2
 =1-2\Rea(tas)+a^2|t|^2(|s|^2+h^2V)
 \le1-\gamma(t),
\]
where
\begin{equation}
 \gamma(t)=2\Rea(tas)-a^2|t|^2
          =|s|^2-a^2\left|t-\frac{\overline s}{a}\right|^2.
 \label{app:decay-disk}
\end{equation}
Completing the square in $t$ explains the conjugate in $C=\overline s/a$. The function $\gamma$ attains its maximum at $C$ and is nonnegative on the closed disk with centre $C$ and radius $|s|/a$. The origin lies on the boundary of this disk, and $\gamma(1)=\alpha\ge0$ puts the endpoint in the disk as well. The path $0\to C\to1$ therefore passes through the point where the bound on the mean squared modulus is smallest.

\begin{figure}[!htb]
\centering
\resizebox{0.80\textwidth}{!}{%
\begin{tikzpicture}[x=30mm,y=30mm,line cap=round,line join=round,>=latex,
                    every node/.style={inner sep=1pt,font=\small}]
  \definecolor{contourfill}{HTML}{ECF3F7}
  \definecolor{contourrim}{HTML}{7799AD}
  \definecolor{contouraxis}{HTML}{869AA7}
  \definecolor{contourfirst}{HTML}{BE552C}
  \definecolor{contoursecond}{HTML}{365F83}
  \definecolor{contourpoint}{HTML}{202D36}
  \path[use as bounding box] (-0.98,-2.25) rectangle (2.57,0.97);
  \coordinate (O) at (0,0);
  \coordinate (C) at ({13/12},{-1/2});
  \coordinate (E) at (1,0);
  \draw[fill=contourfill,draw=contourrim,line width=0.9pt]
    (C) circle[radius={sqrt(205)/12}];
  \draw[contouraxis,thin] (-0.33,0) -- (2.48,0);
  \draw[contouraxis,thin] (0,-1.86) -- (0,0.88);
  \foreach \value in {0,1,2}
    \node[text=contourpoint] at (\value,-1.96) {$\value$};
  \foreach \value in {-1,0}
    \node[text=contourpoint,anchor=east] at (-0.36,\value) {$\value$};
  \node at (1.08,-2.10) {$\operatorname{Re}t$};
  \node at (-0.78,-0.41) {$\operatorname{Im}t$};
  \draw[contourfirst,line width=1.55pt] (O) -- coordinate[pos=0.56] (firstArrow) (C);
  \draw[contoursecond,line width=1.55pt] (C) -- coordinate[pos=0.52] (secondArrow) (E);
  \draw[contourfirst,line width=1.65pt,->] (O) -- (firstArrow);
  \draw[contoursecond,line width=1.65pt,->] (C) -- (secondArrow);
  \foreach \point in {O,C,E}
    \fill[contourpoint] (\point) circle[radius=1.7pt];
  \node[above left] at (-0.04,0.055) {$0$};
  \node[above left] at (0.995,0.05) {$1$};
  \node[right] at (1.20,-0.53) {$C=\overline{s}/a$};
  \node[text=contourfirst] at (0.66,-0.12) {$t=uC$};
  \node[text=contoursecond,anchor=west] at (1.16,0.14) {$t=C+u(1-C)$};
  \node[text=contoursecond] at (1.27,-1.12)
    {$\gamma(t)=|s|^2-a^2|t-C|^2\ge 0$};
\end{tikzpicture}%
}
\caption{The disk $\gamma(t)\ge0$ and the two segments of the contour for $a=3/5$ and $s=13/20+3i/10$. The first segment joins the origin to $C=\overline{s}/a$, where $\gamma$ is largest. The second joins $C$ to the endpoint $1$. The parameter values illustrate the geometry of the estimate.}\label{fig:contour}
\end{figure}
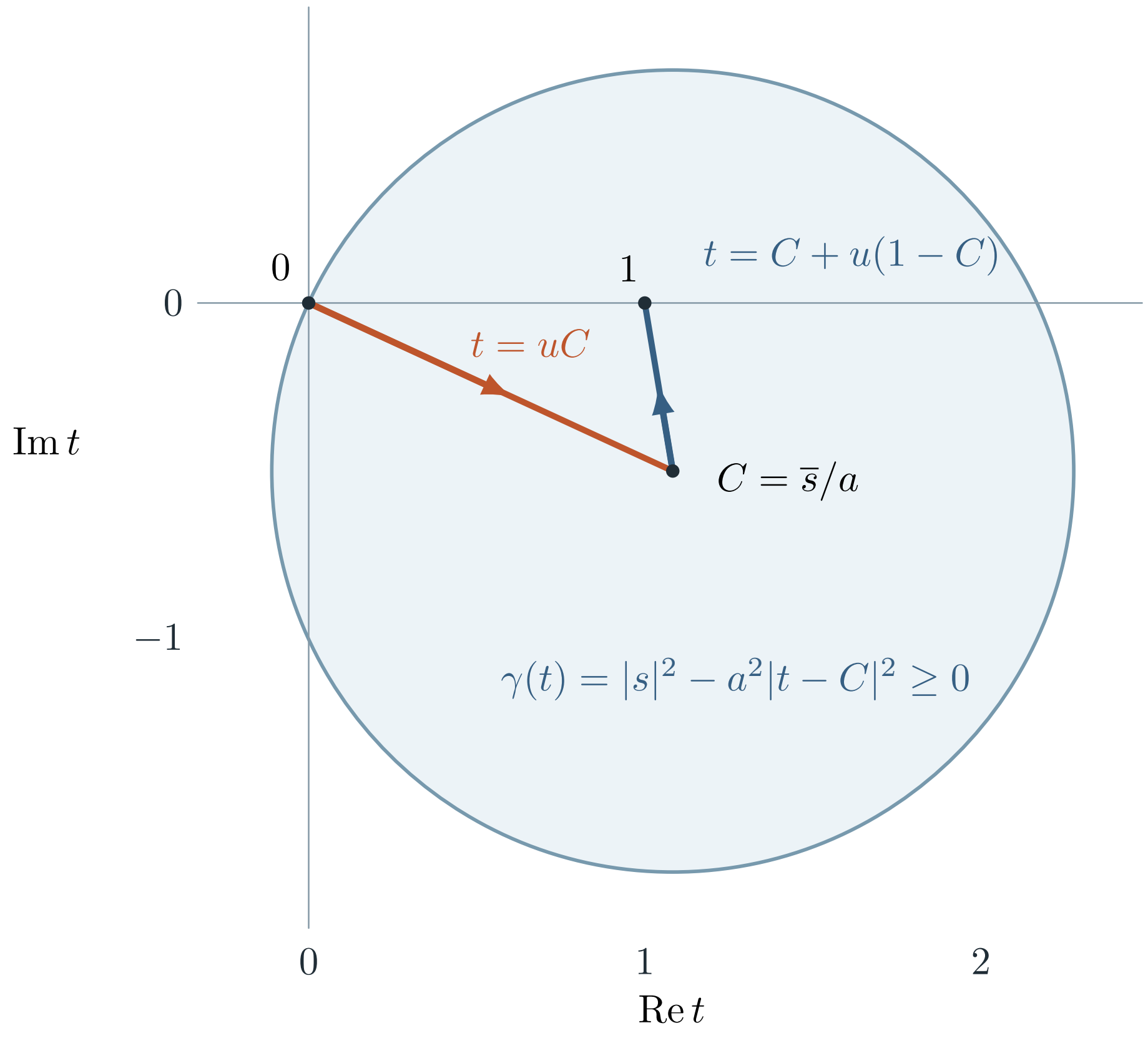

As the integrands are polynomials in $t$, their integrals depend only on the endpoints. We may therefore pass from the real interval to the two segments through $C$ and apply the decay estimate separately on each, whatever the arguments of the individual factors.

On the first segment write $t=uC$, with $0\le u\le1$. Then~\eqref{app:decay-disk} gives
\[
 \gamma(uC)=|s|^2(2u-u^2).
\]
The cancellation in (V2) controls the difference between the product formed with $y_j$ and the one formed with $s$. The factor $t^2$ makes the error small near the origin, where the decay is weakest. Farther along the segment, the arithmetic--geometric mean inequality bounds the product of the remaining $N-2$ factors by a power of their mean squared modulus. This gives the integral involving $u^2e^{-\lambda(u-u^2/2)}$, with $\lambda=N|s|^2$.

On the second segment, put $d=|s|^2-\alpha=a^2|1-C|^2$. For $t=C+u(1-C)$,
\[
 \gamma(t)=\alpha+d(1-u^2)\ge\alpha+d(1-u).
\]
Here both products admit direct decay estimates, and the triangle inequality bounds their difference. The segment has length $\sqrt d/a$ and reduces to a point when $d=0$. Combining its length with the decay along it gives (U3), whose bound decreases exponentially in $v=N\alpha$ up to the factor displayed there.

When $v$ is large, these estimates and the bound at the origin force $a|s|$ close to one and the variance close to zero. The first coefficient of the polynomial then controls the argument of $s$ and places the configuration in the disk $|1-as|\le1/64$. The near-one result identifies the boundary equality case, giving $v<20$ for an interior marked root. For $n\ge257$, reflection forces $v>2\log(N/2)$ and, in the remaining range $v<20$, also gives $a<4/5$. Section~\ref{sec:large} uses these two restrictions in the contour estimate at the origin to obtain the contradiction. For the remaining degrees we use the same two integrals and keep the powers that occur in their bounds.

\section{The finite range of degrees}\label{app:finite}
One polynomial moment formula suffices for the reflected integral and for both segments of the contour. In the reflection argument the exponent is rounded \emph{up}. Here we need only determine whether the integral reaches one, so monotonicity of the power means on a probability space applies even when the base exceeds one. Rounding down pointwise would lose the implication needed for this test. In the two integrals obtained from the evaluation at the origin, the bases lie in the unit interval. Rounding their exponents \emph{down} therefore gives the required upper bounds for the integrals themselves.

The initial domain consists of an integer interval of degrees and a rectangle in $(a,m)$. Each exclusion applies throughout one of the resulting regions, for every integer degree in its interval and every point in its rectangle. For example, the test at the origin succeeds whenever the lower bound for $1-am$ exceeds the sum of the upper bounds for $q^n,E_0,E_1$ throughout that region. At every subdivision the children cover their parent. It follows by induction that the inequalities at the terminal regions exclude the entire initial domain. The accompanying programs reproduce the subdivision and its counts in exact arithmetic. The \hyperref[verification]{discussion of the formal proof} describes how Lean checks these same bounds and uses them in the proof of Theorem~\ref{main}.

\section*{Acknowledgements}
I am grateful to my teacher, \href{https://scholar.google.com/citations?user=HGrSa8QAAAAJ&hl=en}{Professor Petko Hr.~Petkov}, for what he taught me about robust control, singular values, and advanced linear algebra. It was through his teaching that I came to see Sendov's question in terms of the stability of a discrete-time system.

\end{document}